\documentclass[a4paper, 11pt]{article}

\usepackage[utf8]{inputenc} % Encodage des caractères
\usepackage[T1]{fontenc}    % Encodage des fontes
\usepackage[english]{babel} % Langue du document
\usepackage{amsmath, amssymb, amsthm} % Packages AMS pour les mathématiques
\usepackage{geometry}       % Pour régler les marges
\usepackage{graphicx}       % Pour inclure des images
\usepackage{hyperref}       % Pour les liens hypertextes et les références
\usepackage{tocloft}        % Pour personnaliser le sommaire
\usepackage{braket}
\usepackage{cleveref}
\usepackage{tikz}
\usepackage{nicematrix}
\usepackage[numbers]{natbib}

\newcommand{\bk}[2]{\braket{#1|#2}}
\newcommand{\BK}[3]{\braket{#1|#2|#3}}

\numberwithin{equation}{section}

\title{Interlacing Eigenvectors of Large Gaussian Matrices}
\author{Elie Attal, Romain Allez\thanks{Qube Research \& Technologies, Email: elie.attal@polytechnique.edu, romain.allez@qube-rt.com}}
\date{\today}

\begin{document}

% Page de titre
\maketitle

\begin{abstract}
We consider the eigenvectors of the principal minor of dimension $n< N$ of the Dyson Brownian motion in $\mathbb{R}^{N}$ and investigate their asymptotic overlaps with the eigenvectors of the full matrix in the limit of large dimension. We explicitly compute the limiting rescaled mean squared overlaps in the large $n\,, N$ limit with $n\,/\,N$ tending to a fixed ratio $q\,$, for any initial symmetric matrix $A\,$. This is accomplished using a Burgers-type evolution equation for a specific resolvent. In the GOE case, our formula simplifies, and we identify an eigenvector analogue of the well-known interlacing of eigenvalues. We investigate in particular the case where $A$ has isolated eigenvalues. Our method is based on analysing the eigenvector flow under the Dyson Brownian motion.
\end{abstract}

%\tableofcontents

\section{Introduction}
\label{sec:Intro}
\citestyle{numbers}

\indent

Suppose $X \in \mathcal{M}_N(\mathbb{R})$ is symmetric. The Cauchy Interlace Theorem (\cite{haemers1995interlacing}) states that, for any $n < N\,$, the eigenvalues $\mu_1 \geq ... \geq  \mu_n$ of $\tilde{X}$, the principal submatrix of size $n\,$, interlace the eigenvalues $\lambda_1 \geq ... \geq \lambda_N$ of $X$, in the sense that for any $1 \leq i \leq n\,$, we have $\lambda_{i + N - n} \leq \mu_i \leq \lambda_i\,$. We are interested in the relationship between the eigenvectors $\Psi_1\,,...\,, \Psi_N$ of $X$ and the eigenvectors $\Phi_1\,,...\,, \Phi_n$ of $\tilde{X}$. To the best of our knowledge, there is no general result in the literature on this relationship\footnote{Recently, an \textit{eigenvalue-eigenvector identity} has been rediscovered, linking the eigenvectors of a principal submatrix (obtained by deleting one row and one column) to the eigenvalues of the full matrix (see \cite{denton2019eigenvectors}). However, it does not link the eigenvectors of the submatrix to those of the full matrix.}, extending the Cauchy Interlace Theorem to eigenvectors, for general deterministic matrices. We investigate here this relationship for large random matrices with independent entries. More specifically, we consider the case of large random Gaussian matrices.

Let $A$ be an $N \times N$ real symmetric deterministic matrix, and $H_t$ a symmetric random matrix whose coefficients are $H_t^{ii} = \sqrt{\frac{2}{N}} \, B_t^{ii}$ on the diagonal and $H_t^{ij} = H_t^{ji} = \frac{1}{\sqrt{N}} \, B_t^{ij}$ for off-diagonal elements, where $\{B_t^{ij}, \, 1 \leq i \leq j \leq N \}$ are independent Brownian motions. We consider the noisy observation of $A\,$,
\begin{equation*}
X_t := A + H_t \, .
\end{equation*}
For $ n < N$, we define the $N \times N$ matrix $\tilde{X}_t$ by
\begin{equation}
\label{eq:truncation_definition}
\tilde{X}_t^{ij} = 
\begin{cases}
X_t^{ij}\,, \,\,\, \textit{if $i \leq n$ and $j \leq n\,$,} \\
0\hspace{0.48cm}, \,\,\, \textit{otherwise.}
\end{cases}
\end{equation}
It corresponds to the principal $n \times n$ minor of $X_t\,$, where we set all other coefficients to zero to ensure both matrices have the same size:
\begin{equation*}
\tilde{X}_t = \begin{pmatrix}
X_t^{11} & \cdots & X_t^{1n} & 0 & \cdots & 0 \\
\vdots & \ddots & \vdots & \vdots & \ddots & \vdots \\
X_t^{1n} & \cdots & X_t^{nn} & 0 & \cdots & 0 \\
0 & \cdots & 0 & 0 & \cdots & 0  \\
\vdots & \ddots & \vdots & \vdots & \ddots & \vdots \\
0 & \cdots & 0 & 0 & \cdots & 0
\end{pmatrix}\,.
\end{equation*}

We assume that $n$ and $N$ are of the same order in the sense that $n\,/\,N \to q$ as $n\,,\,N \to \infty\,$, where $q$ is a real number in $[0\, , 1]\,$.

The eigenvalues of $X_t$ are denoted by $\lambda_1^t
 \geq ... \geq \lambda_N^t$ and the associated unit norm eigenvectors by $\Psi_i^t\,$, for any $1 \leq i \leq N\,$.
For $\tilde{X}_t\,$, we shall use the notations $\mu_1^t \geq ... \geq \mu_n^t$ for its non-zero eigenvalues, since $\mu_{n+1}^t = ... = \mu_N^t = 0$ almost surely, and $\Phi_i^t$ for its eigenvectors, with $1 \leq i \leq N\,$.

We are interested in studying the overlaps between the eigenvectors of $\tilde{X}_t$ (associated with non-zero eigenvalues) and those of $X_t\,$, i.e. $\bk{\Phi_i^t}{\Psi_j^t}\,$ for $1 \leq i \leq n$ and $1 \leq j \leq N\,$..

For a fixed $t\,$, $H_t$ belongs to the Gaussian Orthogonal Ensemble (GOE), a central model in Random Matrix Theory. The joint probability distribution of the eigenvalues is well-known, and the spectral density converges to the Wigner semicircle as $N \to \infty$ (\cite{anderson2010introduction, potters2020first, tao2012topics}). Additionally, the eigenvectors are uniformly distributed on the unit sphere in $\mathbb{R}^N$, and the rotational invariance of this model implies that they are independent of the eigenvalues (\cite{anderson2010introduction, bai2010spectral, o2016eigenvectors}).

In our model, we introduce time dependence as $H_t$'s coefficients are rescaled standard Brownian motions, making $(X_t)_{t \geq 0}$ a diffusion process. The time parameter $t$ can be seen as a way to select a specific variance for the Gaussian noise. However, it also serves as a powerful way to obtain evolution equations for the objects under investigation. Indeed, in 1962, Dyson demonstrated that the eigenvalues in such a model evolve according to a diffusion known as the Dyson Brownian motion (\cite{dyson1962brownian, potters2020first, tao2012topics, grabiner1999brownian}). It mirrors the dynamics of $N$ particles interacting via a two-dimensional Coulomb potential, characterized by a repulsion force inversely proportional to the particles' distances, and subject to thermal noise. It is a powerful and elegant tool, as it has been used to show universality properties of the eigenvalues of generalized Wigner matrices (\cite{bourgade2016fixed}), as well as to study chaotic billards and disordered metals (\cite{beenakker1994random}). This matrix process also has an associated eigenvector process, which will be the main focus of our study. It has allowed proving universality of the distribution of the eigenvectors of generalized Wigner matrices (\cite{bourgade2017eigenvector}) and has been used in the article that inspired our method to compute the overlaps between the eigenvectors of $X_t$ and those of the initial matrix $A$ (\cite{allez2014eigenvectors}).

The type of noise we use, i.e. additive, is very common in Random Matrix Theory as a model for studying how noise impacts various objects such as Hamiltonians in Quantum Mechanics (\cite{deutsch1991quantum, dyson1962brownian}), Hessian matrices in Machine Learning (\cite{granziol2022learning}), covariance matrices (see Section 5 of \cite{allez2014eigenvectors}), or weight matrices of graphs (\cite{bolla2008noisy}). Therefore, studying the overlaps between $X_t$'s eigenvectors and those of one of its principal minors helps us understand how noise affects the interaction between the entire system and a subsystem in these contexts. For instance, our results allow us to study the typical overlaps between the principal components analysis (PCA) of the covariance matrix of market returns and the PCA of a subset of assets, in cases where there is no real underlying correlation structure. Moreover, in Quantum Mechanics, classical numerical methods involve selecting a finite part of an infinite Hamiltonian for simulation (see Section 8 of \cite{gatti2017applications}). We believe our results help in understanding how noise impacts the information contained in each eigenvector of this submatrix for large dimensions. Finally, minors of Random Matrices have attracted increasing attention in recent years due to their applications in various fields such as compressed sensing (\cite{candes2005decoding, cai2021asymptotic, hu2023extreme, feng2022principal}), percolation theory (\cite{adler2013random}), queueing models (\cite{johansson2006eigenvalues, baryshnikov2001gues}) and conditional independence tests in covariance matrices (\cite{drton2008moments}). Moreover, the eigenvalue interlacing between minors of Wigner matrices has been studied in the microscopic regime ($n = N-k$ for fixed $k$) and shown to converge to a Markov process (see \cite{gorin2017interacting, najnudel2021bead, huang2022eigenvalues}). While these studies primarily focus on the asymptotic properties of the eigenvalues of minors, considerably less attention has been given to their eigenvectors (see, for example, \cite{najnudel2021eigenvector}). Notably, in studying the asymptotics of eigenvalue statistics of submatrices of Wigner matrices, the authors of \cite{li2015central} derived the limit of a statistic generalizing (\ref{eq:S_definition}), using a different method that expresses this limit with an infinite sum of Chebyshev polynomials. We believe that our study can provide new insights into these fields by quantifying the information contained in the eigenvectors of a minor of a noisy matrix, utilizing a novel approach based on the eigenvectors' Itô dynamics. 

Our method is based on the tools developed in \cite{allez2014eigenvectors}: we utilize the Dyson Brownian motion and its associated eigenvector flow to derive evolution equations for the Stieltjes transforms of the objects under investigation. Fortunately, those equations can be solved analytically in the scaling limit, which allows us to obtain closed-form formulas for the limiting mean squared overlaps. More precisely, in Section 2, we introduce the Dyson dynamics of $X_t$ and $\tilde{X}_t$'s eigenvalues and eigenvectors. In our case, the different Brownian motions present have a specific correlation structure that we are able to derive. Moreover, we define the limiting eigenvalue densities along with their Stieltjes transforms and specify the scaling limit that underlies our work. Section 3 is the main part of this paper, where we detail our method for eigenvectors in the bulk of both spectra. We introduce the random function of the complex variables $z$ and $\tilde{z}\,$:
\begin{equation*}
S^{(N)}(z, \tilde{z}, t) := \frac{1}{N} \, \sum_{i=1}^n \, \sum_{j=1}^N \,  \frac{\bk{\Phi_i^t}{\Psi_j^t}^2}{(\tilde{z} - \mu_i^t)(z - \lambda_j^t)} \, .
\end{equation*}
We are able to demonstrate that it is a self-averaging quantity as it converges to a deterministic function that satisfies a specific differential equation (\ref{eq:S_equation}). After solving it using the method of characteristics and inverting the solution with a Stieltjes inversion formula, we obtain the explicit formula (\ref{eq:general_solution}) for the limiting rescaled mean squared overlaps in the case of a general initial matrix $A\,$. In the case where we only observe centered Gaussian noise (i.e. $A \equiv 0$), our formula simplifies and gives the limiting rescaled mean squared overlaps at time $t$ as (see (\ref{eq:solution_noise_only}))
\begin{equation*}
N \, \mathbb{E}\left[\bk{\Phi_{i_n}^t}{\Psi_{j_N}^t}^2\right] \, \longrightarrow \, \frac{(1 - q)\,t}{(1-q)^2 \,t + (\mu - \lambda)\,(\mu - q \lambda)} \,,
\end{equation*}
as $n\,,\,N \to \infty$ with $n\,/\,N \to q\,$, as well as $\lambda_{j_N}^t \to \lambda$ and $\mu_{i_n}^t \to \mu\,$. Using this explicit form, we show that the eigenvectors asymptotically exhibit an interlacing property, in a sense specified in \ref{subsec:Noisy_case}, respecting Cauchy's bounds. Section 4 is an application of our method to the case of a spectrum containing isolated eigenvalues, for which we are also able to obtain the closed-form formulas (\ref{eq:f_equation_simplified}) and (\ref{eq:g_solution}) for the limiting overlaps. Finally, Section 5 extends our work to the case of large random matrices with Bernoulli coefficients. This case is slightly different from the Gaussian one as the isolated eigenvalues diverge in the scaling limit. We start by noting our bulk formula is still numerically valid in that context, before demonstrating how our results can be adapted to obtain a rate of convergence for the overlap of the isolated eigenvectors.

\paragraph{Acknowledgements} We thank Dominik Schröder for his advice and for pointing us to articles \cite{li2015central, huang2022eigenvalues}.
We are also grateful to Adrien Hardy, Benjamin De Bruyne and Enzo Miller for insightful discussions.

\section{Dyson Brownian motion dynamics}
\label{sec:Dyson}
In our model, the eigenvalue and eigenvector processes of $X_t$ are described by the following Dyson Brownian motions:
\begin{equation*}
d\lambda_j^t = \sqrt{\frac{2}{N}} \, dB_j(t) + \frac{1}{N}\,\sum_{\substack{k=1 \\ k \neq j}}^N \, \frac{1}{\lambda_j^t - \lambda_k^t}\,dt \, ,
\end{equation*}
\begin{equation*}
d\Psi_j^t = -\frac{1}{2N} \, \sum_{\substack{k=1 \\ k \neq j}}^N \, \frac{1}{(\lambda_j^t - \lambda_k^t)^2} \, \Psi_j^t \, dt + \frac{1}{\sqrt{N}} \, \sum_{\substack{k=1 \\k \neq j}}^N \, \frac{dW_{jk}(t)}{\lambda_j^t - \lambda_k^t} \, \Psi_k^t \, . 
\end{equation*}
This holds for any $1 \leq j \leq N\,$. Here, $\{B_j\, , \, 1 \leq j \leq N \}$ and $\{W_{jk}\, , \, 1 \leq j \neq k \leq N \}$ are two families of independent (up to a symmetry for $W$) standard Brownian motions, independent of each another. Proofs of these dynamics can be found in Dyson's original paper \cite{dyson1962brownian}, as well as in \cite{tao2012topics} and \cite{potters2020first}, the latter of which uses an interesting perturbation approach.

It is worth noting that the independence of $B$ and $W$ allows us to view the eigenvector dynamics of $X_t$ as a diffusion process in a random environment determined by its eigenvalue trajectories.

For $\tilde{X}_t\,$, note that it can be expressed as $\tilde{A} + \tilde{H}_t$ where $\tilde{A}$ (respectively $\tilde{H}_t$) is the $n$-truncated version of $A$ (respectively $H_t$) filled with zeros, as defined for $\tilde{X}_t$ in (\ref{eq:truncation_definition}). Thus, its eigenvectors have the form $( \tilde{\Phi}^T \, 0\,...\,0)^T$ with $N-n$ zeros and $\tilde{\Phi}$ in $\mathbb{R}^n$ following the same type of dynamics as the eigenvectors of $X_t\,$. Consequently, the dynamics are identical among the eigenvectors associated with non-zero eigenvalues, without interaction with the null space: 
\begin{equation*}
d\mu_i^t = \sqrt{\frac{2}{N}}\,d\tilde{B}_i(t) + \frac{1}{N}\,\sum_{\substack{l= 1 \\ l \neq i}}^n\, \frac{1}{\mu_i^t - \mu_l^t}\,dt \, ,
\end{equation*}
\begin{equation*}
d\Phi_i^t = -\frac{1}{2N \,} \sum_{\substack{l = 1 \\ l \neq i}}^n \, \frac{1}{(\mu_i^t - \mu_l^t)^2}\, \Phi_i^t \, dt + \frac{1}{\sqrt{N}}\, \sum_{\substack{l = 1 \\ l \neq i}}^n \,\frac{d\tilde{W}_{il}(t)}{\mu_i^t - \mu_l^t} \, \Phi_l^t \, , 
\end{equation*}
for any $1 \leq i \leq n\,$. Here, $\tilde{B}$ and $\tilde{W}$ are also independent of each other.

The study of the overlaps $\bk{\Phi_i^t}{\Psi_j^t}$ will require mixing the Brownian motions $W$ and $\tilde{W}$, therefore it is pivotal to derive their correlation. Applying Itô's lemma to the identity $X_t  \, \Psi_j^t = \lambda_j^t \, \Psi_j^t$ gives
\begin{equation}
\label{eq:ito_correlation}
dX_t \,  \Psi_j^t + X_t \, d\Psi_j^t + dX_t \, d\Psi_j^t = d\lambda_j^t \, \Psi_j^t + \lambda_j^t \, d\Psi_j^t \, ,
\end{equation}
where we used the fact that $d\lambda_j^t \, d\Psi_j^t = 0$ due to the independence between $B$ and $W\,$.
Projecting this equation onto $\Psi_k^t$ with $k \neq j$ leads to
\begin{equation}
\label{eq:dW_alternative_form}
\frac{1}{\sqrt{N}} \, dW_{jk}(t) = \bk{\Psi_k^t}{dX_t  \,\Psi_j^t} + \bk{\Psi_k^t}{dX_t \, d\Psi_j^t} \, .
\end{equation}
Similarly, for any $1 \leq i \neq l \leq n\,$,
\begin{equation}
\label{eq:dW_tilde_alternative_form}
\frac{1}{\sqrt{N}} \, d\tilde{W}_{il}(t) = \bk{\Phi_l^t}{d\tilde{X}_t \, \Phi_i^t} + \bk{\Phi_l^t}{d\tilde{X}_t \, d\Phi_i^t} \, ,
\end{equation}
which gives
\begin{equation*}
dW_{jk}(t) \, d\tilde{W}_{il}(t) = N \bk{\Psi_k^t}{dX_t \, \Psi_j^t} \bk{\Phi_l^t}{d\tilde{X}_t \, \Phi_i^t} \, .
\end{equation*}
Finally, a straightforward calculation detailed in Appendix \ref{subsec:correlation} shows that for any $1 \leq j\,,k \leq N$ and $1 \leq i\,,l \leq n\,$, we have
$$
\bk{\Psi_k^t}{dX_t \, \Psi_j^t} \bk{\Phi_l^t}{d\tilde{X}_t \, \Phi_i^t} = \frac{1}{N} \left( \bk{\Phi_i^t}{\Psi_j^t} \bk{\Phi_l^t}{\Psi_k^t} + \bk{\Phi_i^t}{\Psi_k^t} \bk{\Phi_l^t}{\Psi_j^t} \right) dt\,.
$$
Thus, we end up with the correlation
\begin{equation}
\label{eq:correlation}
dW_{jk}(t) \, d\tilde{W}_{il}(t) = \left( \bk{\Phi_i^t}{\Psi_j^t} \bk{\Phi_l^t}{\Psi_k^t} + \bk{\Phi_i^t}{\Psi_k^t}  \bk{\Phi_l^t}{\Psi_j^t} \right) dt \, .
\end{equation}
Since our work focuses on the eigenvectors, this correlation is the most significant. However, in Appendix \ref{subsec:correlation}, we also compute the correlations $dB \, d\tilde{B}\,$, $dB \, d\tilde{W}$ and $dW \, d\tilde{B}\,$. In particular, these correlations are not zero. This is an important fact because, in \cite{allez2014eigenvectors}, the authors take advantage of the independence between $B$ and $W$ to work conditionally on the eigenvalue trajectories. Since $W$ (respectively $\tilde{W}$) is not independent of $\tilde{B}$ (respectively $B$), it would not be possible in our case.

We will first focus on the case of eigenvectors associated with eigenvalues in the bulk of the spectra, which requires us to assume two properties in our model:
\begin{enumerate}
\item $A := (A_N)_{N \geq 0}$ is a sequence of deterministic $N \times N$ real symmetric matrices, with spectrum $\lambda_i^0 \geq ... \geq \lambda_N^0$\,. We assume that its empirical eigenvalue density converges weakly towards a continuous probability distribution on $\mathbb{R}\,$, denoted by $\rho(\cdot,0)\,$, in the sense that
\begin{equation*}
\frac{1}{N} \,\sum_{j = 1}^N\, \delta_{\lambda_j^0}(d\lambda) \,\longrightarrow \, \rho(\lambda, 0)\,d\lambda \,.
\end{equation*}
\item Similarly, we assume for the non-zero part of $\tilde{A}$'s spectrum:
\begin{equation*}
\frac{1}{n} \, \sum_{i = 1}^n \, \delta_{\mu_i^0}(d\mu) \,\longrightarrow\, \tilde{\rho}(\mu, 0) \, d\mu \,.
\end{equation*}
\end{enumerate}
We use the notations $\rho(\lambda,t)$ for the limiting spectral density of $X_t$ and $\tilde{\rho}(\mu,t)$ for the limiting density of the non-zero part of $\tilde{X}_t$'s spectrum. Note that these limits are deterministic, as we expect the eigenvalues to stick to their typical quantile positions in the large $N$ limit. Some natural tools to consider are their Stieltjes transforms, defined for $z \in \mathbb{C} \setminus \mathbb{R}$ by 
\begin{align}
\label{eq:G_definition}
G_N(z,t) := \frac{1}{N} \,\sum_{j=1}^N \, \frac{1}{z - \lambda_j^t} \, &\longrightarrow \, G(z,t) := \int_{\mathbb{R}} \frac{\rho(\lambda,t)}{z - \lambda}\,d\lambda \,, \\
\label{eq:G_tilde_definition}
\tilde{G}_n(z,t) := \frac{1}{n} \,\sum_{i = 1}^n \, \frac{1}{z - \mu_i^t} \, &\longrightarrow \, \tilde{G}(z,t) := \int_{\mathbb{R}} \frac{\tilde{\rho}(\mu,t)}{z - \mu}\, d\mu \,.
\end{align}
These transforms converge to deterministic limits as they are typically self-averaging objects (see \cite{potters2020first, tao2012topics}). They are convenient tools for studying the associated spectral densities. By projecting onto the real axis, we can then reconstruct the density using the Sokhotski-Plemelj formula
\begin{equation}
\label{eq:Stieltjes_inversion}
\lim_{\varepsilon \to 0^+} \, G(\lambda \pm i \,\varepsilon, t) = v(\lambda,t)  \, \mp \, i \pi \rho(\lambda,t) \,,
\end{equation}
where $v(\lambda,t) := P.V. \int \frac{\rho(\lambda',t)}{\lambda - \lambda'}\,d\lambda'$ is the Hilbert transform of $\rho$ and $P.V.$ denotes Cauchy's principal value.

By applying Itô's lemma, \cite{voiculescu1986addition, allez2014eigenvectors, potters2020first} and \cite{tao2012topics} find the classical Burgers evolution equation in the scaling limit,
\begin{equation}
\label{eq:G_burgers}
\partial_t G = - G \, \partial_z G \,.
\end{equation}
In our case, a similar equation for $\tilde{G}$ is easily obtained:
\begin{equation}
\label{eq:G_tilde_burgers}
\partial_t  \tilde{G} = -q  \tilde{G} \, \partial_z \tilde{G} \, .
\end{equation}
In Appendix \ref{subsec:burgers}, we demonstrate how, using the method of characteristics, one can derive from these Burgers equations the following implicit equations depending on the initial conditions,
\begin{align}
\label{eq:G_implicit}
G(z,t) &= G \left(z - t\,G(z,t)\,,0 \right) \,, \\
\label{eq:G_tilde_implicit}
\tilde{G}(z,t) &= \tilde{G} \left( z - qt \,\tilde{G}(z,t)\,, 0 \right) \,.
\end{align}
The relation (\ref{eq:G_implicit}) was first derived by Biane in 1997, see Proposition 2 in \cite{biane1997free}\,. These objects and the equations they satisfy will play a central role in our study.

Our goal is to investigate the limiting behaviour of the overlaps between the eigenvectors of $\tilde{X}_t$ and those of $X_t$, i.e. $\bk{\Phi_i^t}{\Psi_j^t}\,$. Recalling that the eigenvectors are defined up to a sign, it makes more sense to consider the square of this quantity, which also satisfies useful normalisation constraints.
Even though in the large $N$ limit the eigenvalue density becomes deterministic (see \cite{potters2020first, tao2012topics, bai2010spectral} for convergence towards Wigner's semicircle in the GOE case and \cite{bourgade2016fixed} for results on generalized Wigner matrices), the eigenvectors and their projections are still random objects in the scaling limit in the sense that they converge in distribution towards Gaussian variables in the GOE case (\cite{tao2012random, o2016eigenvectors}), and in the generalized Wigner case (\cite{bourgade2017eigenvector}). Therefore, we will consider the expectation $\mathbb{E}[\cdot]$ of the squared overlaps, over the whole randomness of $B\,$, $W\,$, $\tilde{W}$ and $W\,$, i.e. $\mathbb{E}\left[ \bk{\Phi_i^t}{\Psi_j^t}^2 \right]$. Finally, the normalization constraint $\sum_{j = 1}^N \bk{\Phi_i^t}{\Psi_j^t}^2 = 1$ indicates this quantity vanishes as $1\,/\,N\,$. Thus, our goal will be to study the limit of
\begin{equation*}
u_{i|j}(t) := N \, \mathbb{E} \left[\bk{\Phi_i^t}{\Psi_j^t}^2\right] \,.
\end{equation*}

The scaling limit for studying the overlaps is as follows: assuming $t$ is independent of $N$ (macroscopic regime), we have $n\,, N \to \infty$ with $n\,/\,N \to q\,$. Additionally, we consider indices of the form $j_N \,/\, N \to y \in [0\,,1]$ for the matrix $X_t$ and $i_n \,/\, n \to x \in [0\,,1]$ for $\tilde{X}_t\,$. This keeps the distances between eigenvalues of order 1 instead of $1 \,/\, N$ if we had used fixed indices. In this limit, we have $\lambda_{j_N}^t \to \lambda(y,t)$ (respectively $\mu_{i_n}^t \to \mu(x,t)$) where $\lambda(\cdot, t)$ (respectively $\mu(\cdot, t)$) is the quantile function of $\rho(\cdot,t)$ (respectively $\tilde{\rho}(\cdot,t)$). These are defined by
\begin{equation}
\label{eq:quantile_function}
x = \int_{\lambda(x,t)}^{+\infty} \rho(\lambda,t)\, d\lambda = \int_{\mu(x,t)}^{+\infty} \tilde{\rho}(\mu,t) \, d\mu \,,\quad \text{for any $x \in [0\,,1]\,$.} 
\end{equation}
In this scaling limit, we expect $u_{i_n|j_N}(t)$ to converge to a function $U(x,y,t)$ that we want to make explicit, in the same way it is done in \cite{allez2014eigenvectors}.

\section{Limiting rescaled mean squared overlaps in the bulk}
\label{sec:Overlaps}
\subsection{The general case}
For readability, we introduce the (non-symmetric) notation $\bk{i}{j} := \bk{\Phi_i^t}{\Psi_j^t}\,$. Applying Itô's lemma and using the correlation (\ref{eq:correlation}) (see Appendix \ref{subsec:App_overlaps}), we obtain
\begin{align}
\label{eq:ito_squared_overlaps}
d\left(\bk{i}{j}^2\right) &= \frac{1}{N} \,\sum_{\substack{k = 1 \\ k \neq j}}^N \, \frac{\bk{i}{k}^2 - \bk{i}{j}^2}{(\lambda_j^t - \lambda_k^t)^2} \,dt + \frac{1}{N} \,\sum_{\substack{l = 1 \\ l \neq i}}^n \, \frac{\bk{l}{j}^2 - \bk{i}{j}^2}{(\mu_i^t - \mu_l^t)^2} \,dt \nonumber \\
& \quad  + \frac{2}{N} \, \sum_{\substack{l = 1 \\ l \neq i}}^n \, \sum_{\substack{k = 1 \\ k \neq j}}^N \, \frac{\left( \bk{i}{j} \bk{l}{k} + \bk{i}{k} \bk{l}{j} \right)^2}{(\mu_i^t - \mu_l^t)(\lambda_j^t - \lambda_{k}^t)}\,dt \nonumber \\
& \quad + \frac{2}{\sqrt{N}}\, \sum_{\substack{k = 1 \\ k \neq j}}^N \, \frac{dW_{jk}(t)}{\lambda_j^t - \lambda_k^t}\bk{i}{k} \bk{i}{j} + \frac{2}{\sqrt{N}} \, \sum_{\substack{l=1 \\ l \neq i}}^n \, \frac{d\tilde{W}_{il}(t)}{\mu_i^t - \mu_l^t}\bk{l}{j} \bk{i}{j} \,.
\end{align}
The first two terms of this equation recall the one obtained in \cite{allez2014eigenvectors}. Note that $\frac{1}{N} \,\sum_{k \neq j}\, \frac{1}{(\lambda_j - \lambda_k)^2}$ diverges in the scaling limit for eigenvalues in the bulk. Therefore, the compensation brought by the factor $\bk{i}{k}^2 - \bk{i}{j}^2$ is essential in order to have a limiting object. We also remark that the factor $(\bk{i}{j} \bk{l}{k} + \bk{i}{k} \bk{l}{j})^2$ poses a challenge: in the large $N$ limit, as the eigenvalues are expected to become deterministic, taking the expectation of (\ref{eq:ito_squared_overlaps}) would not yield an autonomous equation for the function $U$ (unlike \cite{allez2014eigenvectors}). One would have to assume that overlaps on different eigenvectors are asymptotically independent to write $N^2 \, \mathbb{E}\left[\bk{i_n}{j_N}^2 \bk{l_n}{k_N}^2 \right] \to U(x,y,t) \,U(x', y', t)$ where $(i_n / n \,,\,  l_n/n \,,\, j_N /N\,,\, k_N/N) \to (x\,,\,x'\,,\,y\,,\,y')\,$. We address this issue by initially working with the random squared overlaps $\bk{i}{j}^2$ rather than the rescaled mean squared overlaps $u_{i|j}\,$. Let us introduce the complex function
\begin{equation}
\label{eq:S_definition}
S^{(N)}(z,\tilde{z}, t) := \frac{1}{N} \,\sum_{i=1}^n\, \sum_{j = 1}^N \, \frac{\bk{i}{j}^2}{(\tilde{z}- \mu_i^t)(z - \lambda_j^t)} \,.
\end{equation}
It can be viewed as a double Stieltjes transform where each pole is associated with a different overlap. Similar functions have been previously used in \cite{bun2018overlaps, pacco2023overlaps} and \cite{ledoit2011eigenvectors}. However, in \cite{bun2018overlaps} and \cite{pacco2023overlaps}, the authors use its mean, which could not work in our setup due to the issue mentioned above, forcing us to use its random counterpart. By summing the contributions of the different overlaps, the intuition is that this object is self-averaging in the large $N$ limit, meaning it converges to a deterministic function that is its expectation. We are going to show that this intuition is correct, as this function converges almost surely to the solution of a deterministic differential equation that we are able to solve. This a key result as it demonstrates that our method can be applied to a broad range of problems, even without a clean equation for the mean squared overlaps. More specifically, in Appendix \ref{subsec:stieltjes_equation}, we show that $S\,$, the limit of $S^{(N)}\,$, almost surely verifies
\begin{equation}
\label{eq:S_equation}
\partial_t S = -G(z,t) \, \partial_z S -q\tilde{G}(\tilde{z},t) \, \partial_{\tilde{z}} S + S^2\,,
\end{equation}
where $G$ and $\tilde{G}$ are the Stieltjes transforms of the limiting eigenvalues densities introduced in (\ref{eq:G_definition}) and (\ref{eq:G_tilde_definition}). Note that the characteristics curves of this equation are the same as those in (\ref{eq:G_burgers}) and (\ref{eq:G_tilde_burgers}), which allows us to solve it. Hence, the variables $z - t\,G(z,t)$ and $\tilde{z} - qt \,\tilde{G}(\tilde{z},t)$ should appear in the solution. The rest is determined by the Ricatti part $f' = f^2\,$, for which the unique solution is $f(x) = \frac{f(0)}{1 - f(0)\,x}\,$. The final solution, depending on the initial conditions, $G\,$, and $\tilde{G}$ (see Appendix \ref{subsec:stieltjes_equation} for the resolution using the method of characteristics), is
\begin{equation}
\label{eq:S_general_solution}
S(z, \tilde{z}, t) = \frac{S\left(z - t\,G(z,t)\,,\, \tilde{z} - qt\,\tilde{G}(\tilde{z}, t)\,,\, 0\right)}{1 - t\,S\left(z - t\,G(z,t)\,,\, \tilde{z} - qt\,\tilde{G}(\tilde{z}, t)\,,\, 0\right)} \,.
\end{equation}
Since $S^{(N)}$ converges to a deterministic function, we deduce it is indeed self-averaging and that its limit is also the limit of $\mathbb{E}\left[ S^{(N)} \right]\,$. Recalling that we expect the eigenvalues to become deterministic, the average is asymptotically taken only over the overlaps, giving
\begin{equation*}
S(z, \tilde{z},t) = q \int_0^1 \int_0^1 \frac{U(x,y,t)}{(\tilde{z} - \mu(x,t))(z - \lambda(y,t))}\,dx\,dy \,,
\end{equation*}
where $U\,$, $\lambda(\cdot,t)\,$ and $\mu(\cdot, t)$ are defined at the end of Section \ref{sec:Dyson}. The same reasoning was applied in \cite{bun2018overlaps} for a similar function. By introducing the function $W$ defined by $W\left(\mu(x,t)\,,\, \lambda(y,t)\,,\,t\right) = U(x,y,t)\,$, which is equivalent to using the eigenvalues as indices, this can be rewritten as
\begin{equation*}
S(z,\tilde{z},t) =  q \int_{\mathbb{R}} \int_{\mathbb{R}} \frac{W(\mu, \lambda, t) \, \tilde{\rho}(\mu,t) \, \rho(\lambda,t)}{(\tilde{z} - \mu)(z - \lambda)} \,d\mu\, d\lambda \,.
\end{equation*}
This demonstrates why this choice of function is useful: thanks to the inversion formula derived in \cite{bun2018overlaps}, we can reconstruct our goal function by projecting onto the real axis, similarly to (\ref{eq:Stieltjes_inversion}):
\begin{equation}
\label{eq:double_Stieltjes_inversion}
W(\mu, \lambda, t) = \lim_{\varepsilon \to 0^+} \, \frac{\Re\left[ S(\lambda - i \,\varepsilon, \mu + i \,\varepsilon, t) - S(\lambda - i \,\varepsilon, \mu - i \,\varepsilon, t)  \right]}{2 q \pi^2 \rho(\lambda,t) \tilde{\rho}(\mu,t)} \,.
\end{equation}
Therefore, the general solution to our problem is given by
\begin{equation}
\label{eq:general_solution}
W(\mu, \lambda , t) = \frac{1}{2 q \pi^2 \rho(\lambda,t) \tilde{\rho}(\mu,t)} \, \Re \left[\frac{S(y,\tilde{y}^*,0)}{1 - t S(y,\tilde{y}^*,0)} -  \frac{S(y,\tilde{y},0)}{1 - t S(y,\tilde{y},0)} \right] \,,
\end{equation}
with $y := \lambda - t\,v(\lambda,t) - i \,\pi t \,\rho(\lambda,t)$ and $\tilde{y} := \mu - qt \,\tilde{v}(\mu,t) - i\, q\pi t \, \tilde{\rho}(\mu,t)\,$. We used the notation $\tilde{y}^*$ for the complex conjugate of $\tilde{y}\,$.

Even though $S(\cdot\,,\, \cdot\,,\, 0)$ is fully known with $A\,$, it cannot be expressed simply in the general case. Therefore, it is interesting to consider the simple case $A \equiv 0\,$, which corresponds to the observation of a pure noise matrix. Moreover, this case, equivalent to $A = \alpha \,I_N\,$, provides insights into what constitutes a significant overlap in the context of real applications, by comparing it to the overlap structure of a GOE matrix.

\subsection{The GOE case}
\label{subsec:Noisy_case}
In this section, we derive the simplified form of the previous formula in the case $A \equiv 0\,$, meaning we observe a centered Gaussian matrix. The advantage of this particular case is that we have closed-form solutions for (\ref{eq:G_implicit}) and (\ref{eq:G_tilde_implicit}). Indeed, the initial condition on $G$ is $G(z,0) = 1 \,/\,z\,$, which allows us to compute it as one of the roots of a second-order polynomial. Since we need to have $G(z,t) \sim 1 \,/\,z$ as $z\to \infty\,$, only one root is possible and we get
\begin{equation*}
G(z,t) = \frac{z - \sqrt{z^2 - 4t}}{2t} \,.
\end{equation*}
This is the Stieltjes transform of $\rho\,$, the Wigner semicircular density of radius $2 \sqrt{t}\,$:
\begin{equation*}
\rho(\lambda,t) = \frac{\sqrt{(4t - \lambda^2)_+}}{2 \pi t} \,.
\end{equation*}
Similarly, the implicit equation for $\tilde{G}$ can now be solved using the same initial condition, which shows that $\tilde{\rho}$ is the Wigner semicircular density of radius $2 \sqrt{qt}\,$. Their respective Hilbert transforms are $v(\lambda,t) = \lambda \,/\, 2t$ and $\tilde{v}(\mu,t) = \mu \,/\, 2qt$ for eigenvalues inside the bulk. We can note here that in the scaling limit, the density of the non-zero part of $\tilde{X_t}$'s spectrum is not the same as $X_t$'s as it is scaled by a factor of $\sqrt{q}\,$. This discrepancy arises because both matrices are rescaled by the factor $1 \,/\,\sqrt{N}\,$, and not $1\,/\,\sqrt{n}$ for $\tilde{X}_t\,$.

Also, since all the eigenvalues are null at $t = 0\,$, we have $S(z, \tilde{z}, 0) = q \,/\, z \tilde{z}$ which simplifies (\ref{eq:S_general_solution}) into
\begin{equation}
\label{eq:S_solution_noise_only}
S(z, \tilde{z}, t) = \frac{q}{\left(z - t\,G(z,t)\right)\,\left(\tilde{z} - qt\, \tilde{G}\left(\tilde{z}, t\right)\right) - qt} \,.
\end{equation}
Finally, using our inversion formula (\ref{eq:double_Stieltjes_inversion}) along with the explicit forms of the densities and their Hilbert transforms, we can perform some simplifications (see Appendix \ref{subsec:cauchy_simplification}) and arrive to the final form
\begin{equation}
\label{eq:solution_noise_only}
W(\mu, \lambda, t) = \frac{(1 - q)\,t}{(1 - q)^2\, t + (\lambda - \mu)\,(q \lambda - \mu)} \,.
\end{equation}
This formula is the main result of our paper. In the case $A = \alpha\, I_N\,$, the result is identical except that the eigenvalues are shifted by $\alpha\,$.

Let us first remark that $W$ has the form of a Cauchy-like distribution in $\lambda$ or $\mu\,$, which was already observed in \cite{allez2014eigenvectors} when investigating the overlaps between $X_t$'s eigenvectors and $A$'s.
As q increases, this function becomes more peaked and converges towards $\delta_{\lambda\mu}$ when $q = 1\,$, which corresponds to $\tilde{X}_t = X_t\,$. If $q$ decreases, $\tilde{X}_t$ is "further" from $X_t$ and the overlaps are more uniform, the information from $X_t$ gets lost.
When we increase the noise, which corresponds to increasing $t\,$, the behaviour is similar: $W$ becomes flatter as the missing part from $X_t$ in $\tilde{X}_t$ has more variance. 

We want to use this formula to prove an analogue of the Cauchy Interlace Theorem
\begin{equation*}
\lambda_{i+N-n}^t  \leq \mu_i^t \leq \lambda_i^t \quad , \quad \text{for $1 \leq i \leq n$} \,,
\end{equation*}
for eigenvectors in our model. We first note that the asymptotic version of this relation is obtained by taking indices $i_n$ such that $i_n \, / \, n \to x\,$, which gives
\begin{equation}
\lambda(qx + 1 -q, t) \leq \mu(x,t) \leq \lambda(qx,t) \quad , \quad \text{for $x \in [0 \,, 1]$} \,.
\end{equation}
For a fixed $\mu\,$, associated with an $x$ in $[0\,,1] \,$, we look for $\lambda_*$ such that the function $W(\mu, \cdot, t)\,\rho(\cdot, t)$ reaches its maximum. This corresponds to the region of eigenvalues where the eigenvectors contain the most information about the eigenvector $\Phi$ associated with $\mu\,$. Note that multiplying by the eigenvalue distribution $\rho$ is essential because it accounts for the density of eigenvalues near $\lambda_*\,$. Differentiating this function with respect to $\lambda$ shows that $\lambda_*$ must be a root of a cubic polynomial:
\begin{equation*}
q\,\lambda_*^3 - \left((1 + 6q + q^2)\,t + \mu^2 \right) \lambda_* + 4(1+q)\,t\,\mu = 0 \,.
\end{equation*}
In Appendix \ref{subsec:App_cubic}, we demonstrate that for any $-2\sqrt{qt} \leq \mu \leq 2\sqrt{qt}\,$, there is a unique solution $\lambda_*$ to this equation that lies in $[-2\sqrt{t}\,,2\sqrt{t}]\,$, and that it satisfies
\begin{equation*}
\lambda(qx+1-q,t) \leq \lambda_* \leq \lambda(qx,t) \,.
\end{equation*}
This relation serves as an analogue of the Cauchy Interlace Theorem for eigenvectors, as it shows that for one of $\tilde{X}_t$'s eigenvalues $\mu\,$, the region in $X_t$'s spectrum where the eigenvectors locally contain the most information about the eigenvector associated with $\mu$ is within Cauchy's bounds. Moreover, our result is more precise, as one can derive an analytical formula for $\lambda_*$ (as a solution to a cubic equation) and because we can also justify that $\lambda_*$ lies between $\mu$ and $\mu \,/\, \sqrt{q}\,$, which is a tighter interval (see Appendix \ref{subsec:App_cubic}).

Figure \ref{fig:null_A} shows a comparison of our formula for $W(\mu, \lambda,t) \, \rho(\lambda,t)$ with simulated rescaled mean squared overlaps in the case $A \equiv 0\,$, with $q = 0.9$ which implies very peaked curves. The fit with the data appears to be excellent. These curves also demonstrate that the more extreme the eigenvalue associated with the eigenvector for which we plot the overlaps (i.e. close to the semicircle edge), the more peaked the overlaps curve becomes. Additionally, we represent with coloured bands the interlacing intervals given by $\lambda(qx + 1 -q, t)$ and $\lambda(qx,t)\,$. One can see the maximum of each curve is indeed reached within these bounds.

\begin{figure}[!h]
\centering
\includegraphics[scale=0.7]{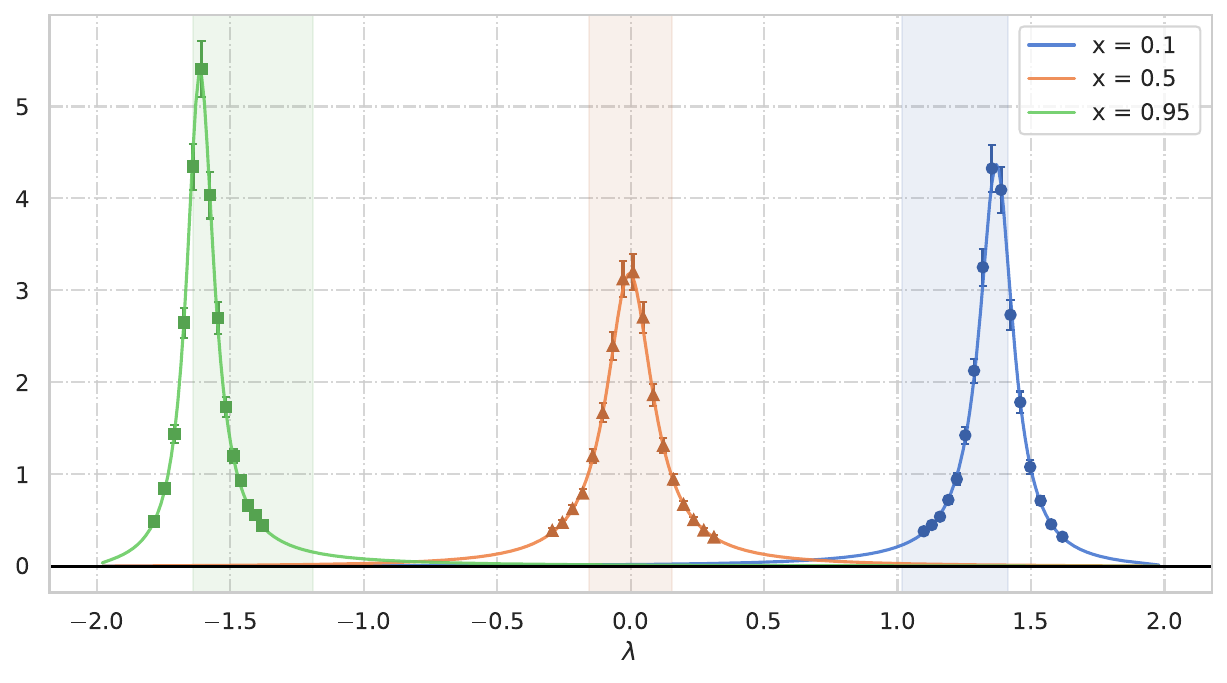}
\caption{Comparison of our theoretical formula for $W(\mu, \lambda,t) \, \rho(\lambda,t)$ with numerical simulations of $N\,\mathbb{E}\left[\bk{\Phi_i^t}{\Psi_j^t}^2\right] \, \rho(\lambda_j^t, t)\,$, as a function of $\lambda$ in the case $N = 500\,$, $t = 1\,$, $q = 0.9$ and $A \equiv 0$ for different values of $\mu$ given by the quantiles $x = 0.1$ (blue circles for data and blue curve for theory), $x = 0.5$ (orange triangles and orange curve) and $x = 0.95$ (green squares and green curve). The data points are shown with 99\% confidence intervals. For each value of $x\,$, the coloured band represents the interlacing interval bounded by $\lambda(qx + 1-q,t)$ and $\lambda(qx,t)\,$.}
\label{fig:null_A}
\end{figure}

\begin{figure}[h!]
	
    \centering
    \begin{minipage}[b]{0.48\textwidth}
        \centering
        \includegraphics[width = \textwidth, height = 5cm]{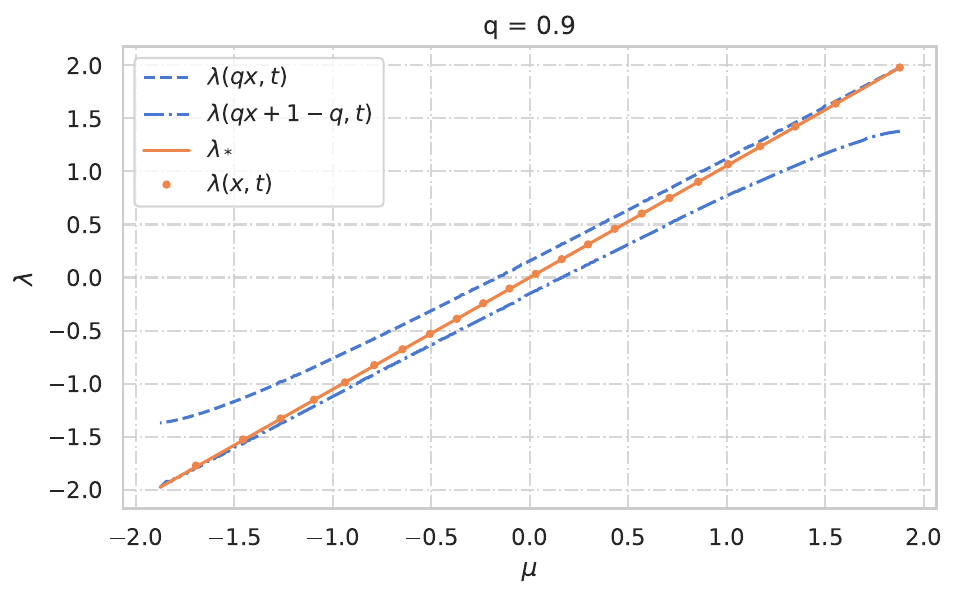}
    \end{minipage}
\hspace{0cm}
    \begin{minipage}[b]{0.48\textwidth}
        \centering
        \includegraphics[width = \textwidth, height = 5cm]{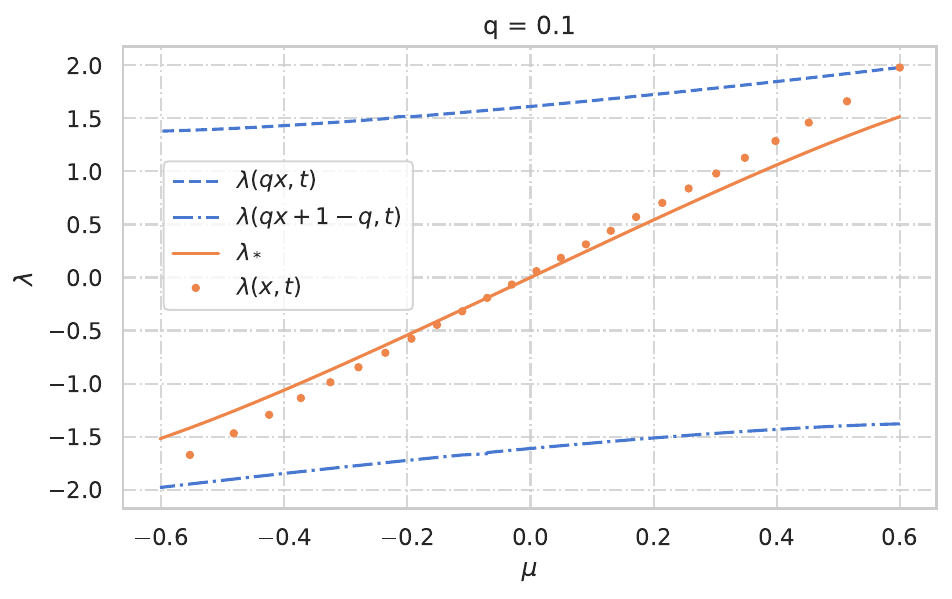}
    \end{minipage}
    \caption{Comparison of $\lambda_*$ (orange curve) with the interlacing bounds (blue curves) as $\mu$ varies. We also plot $\lambda(x,t) = \mu \,/\, \sqrt{q}$ (orange dots) to show it differs from $\lambda_*\,$, which is shifted towards the center of the spectrum. We use the parameters $A = 500$ and $t=1\,$. \textbf{Left:} $q = 0.9\,$. \textbf{Right:} $q = 0.1\,$.}
    \label{fig:Interlacing}
\end{figure}

One could think that asymptotically, $\tilde{X}_t$'s spectrum behaves as a shrunk version of $X_t$'s, meaning that the best projector for the eigenvector associated with $\mu(x,t)$ is the eigenvector of $\lambda(x,t)\,$, positioned at the same quantile in its spectral density. This would imply $\lambda_* = \mu(x,t) \, /\, \sqrt{q}\,$. Although this would still respect Cauchy's interlacing bounds, it is not the case. Indeed, the third-order polynomial for which $\lambda_*$ is a zero differs from zero in $\mu \,/\, \sqrt{q}$ in most cases (except if $\mu = 0$ or $q=1\,$, for instance). In Figure \ref{fig:Interlacing}, we show how $\lambda_*$ varies for different values of $\mu(x,t)\,$. In particular, for $q = 0.1\,$, we clearly see it differs from $\lambda(x,t) = \mu(x,t) \,/\, \sqrt{q}\,$, with a noticeable shift towards the center of the spectrum. This shift is due to the fact that the eigenvalue density $\rho$ is zero at $\pm 2 \sqrt{t}\,$, so there cannot be any concentration of the information of the eigenvector associated with $\mu(x,t)$ at the edge. Moreover, the plot with $q = 0.9$ shows that the bounds become optimal in the limit $q \to 1$ (which can be proved mathematically).

\section{The case of a spiked matrix}
In this section, we extend our analysis to a scenario where the matrix $A$ has a spike, representing an isolated eigenvalue. This situation is common in applications such as finance, where covariance matrices often exhibit spiked behaviour due to significant market factors. The method used is very similar to the one detailed in Section \ref{sec:Overlaps}.

Consider $A = \Psi \Psi^T$, a rank-1 matrix, for a certain vector $\Psi$ in $\mathbb{R}^N$. This matrix has one non-zero eigenvalue $\lambda = \|\Psi\|^2$ associated with the unit norm eigenvector $\pm\, \Psi\, /\, \|\Psi\|\,$. More precisely, $\Psi := (\Psi_N)_{N \geq 0}$ is a sequence of vectors such that $\| \Psi_N \|^2 \to \lambda < \infty\,$. When $A$ is perturbed by the Brownian motion $H_t\,$, the authors of \cite{allez2014eigenvectors} demonstrate that asymptotically, the spectrum will exhibit a bulk following the Wigner semicircle distribution and a spike $\lambda_1(t)$ with a deterministic trajectory. It is obtained by sending $N$ to infinity in the Dyson Brownian motion dynamics:
\begin{equation*}
\frac{d\lambda_1}{dt} = \int_{\mathbb{R}} \frac{\rho(\lambda,t)}{\lambda_1(t) - \lambda} \,d\lambda \,,
\end{equation*}
where $\rho$ corresponds to the spectral density of the bulk, i.e the semicircular density of radius $2 \sqrt{t}\,$. This can be solved to yield
\begin{equation*}
\lambda_1(t) = \lambda + \frac{t}{\lambda} \,.
\end{equation*}
This dynamic holds for any $t < \lambda^2$. Indeed, $t_c := \lambda^2$ is the critical time at which the edge of the spectrum catches the spike. After $t_c\,$, the matrix $X_t$ only has a bulk of eigenvalues. Consequently, we will only consider times $t < \lambda^2$ for the remainder of this section. For finite $N\,$, we will denote $X_t$'s top eigenvalue by $\lambda_1^N(t)\,$.

When considering $\tilde{X}_t\,$, the $n$-truncated version of $X_t = A + H_t\,$, two cases arise:
\begin{itemize}
\item[$\bullet$] \textbf{Spike-spike:} if $\Psi$'s first $n$ coefficients $\Psi^1\,, ...\,, \Psi^n$ are not all zero, then $\tilde{X}_t$ also has a spike $\mu_1^N(t)$ converging to $\mu_1(t) = \mu + q\,\frac{t}{\mu}\,$, where $\mu = \| \Psi^{|n} \|^2$ is the initial eigenvalue of $\tilde{A}$ (we define the $n$-truncated version of an $\mathbb{R}^N$ vector $\Psi$ by $\Psi^{|n} := (\Psi^1 ... \,\Psi^n \, 0\,...\,0)^T$). This holds for any $t < \mu^2 \,/\, q\,$.
\item[$\bullet$] \textbf{Spike-bulk:} if $\Psi^1 = ... = \Psi^n = 0\,$, then $\tilde{A} \equiv 0$ and $\tilde{X}_t$ consists solely of a bulk of eigenvalues, following the Wigner semicircle density of radius $2\sqrt{qt}\,$.
\end{itemize}
In both cases, we want to study how $\Psi_1^t\,$, the eigenvector associated with the spike of $X_t\,$, is projected onto the eigenbasis of $\tilde{X}_t\,$.

From now on, $G(z,t)$ stands for the Stieltjes transform of the Wigner semicircle of radius $2\sqrt{t}\,$:
\begin{equation*}
G(z,t) = \frac{z - \sqrt{z^2 - 4t}}{2t} \,,
\end{equation*}
and $\rho$ (respectively $v$) is the associated density (respectively its Hilbert transform). Similarly, $\tilde{G}\,$, $\tilde{\rho}$ and $\tilde{v}$ are the corresponding functions for the Wigner semicircle of radius $2\sqrt{qt}\,$.
Also, $S$ corresponds to our double Stieltjes transform in the bulk in the case $A \equiv 0\,$, i.e.
\begin{equation*}
S(z, \tilde{z},t) = \frac{q}{(z - t\,G(z,t))\,\left(\tilde{z} - qt \,\tilde{G}\left(\tilde{z},t\right)\right) - qt} \,.
\end{equation*}

\subsection{Spike-spike overlap}
\label{subsec:spike_spike}
\begin{figure}[!h]
\centering
\includegraphics[scale=0.9]{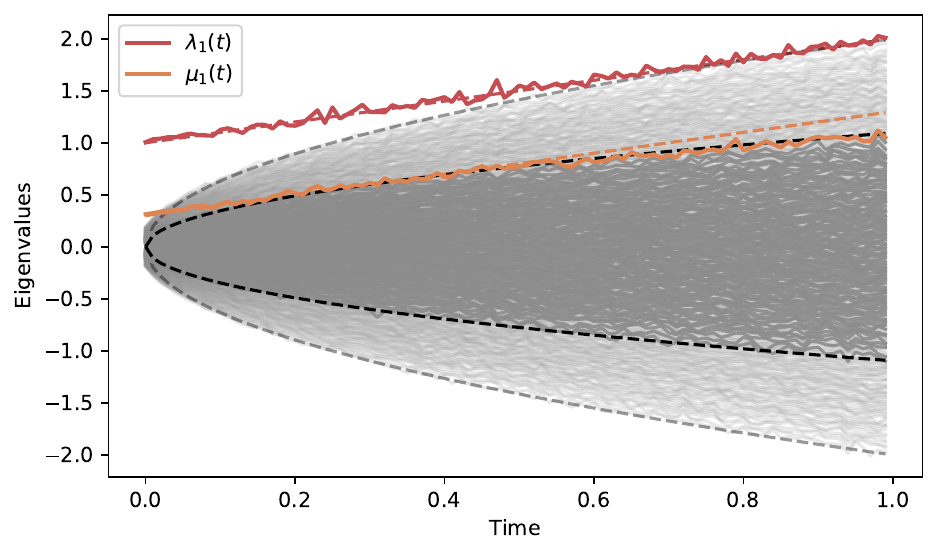}
\caption{Illustration of the time evolution of the two spikes and bulks with $\Psi = (1\,...\,1)^T \,/\, \sqrt{N}$ and $q = 0.3\,$, meaning $\lambda = 1$ and $\mu = 0.3\,$. The simulation was done with $N = 300\,$. We can see that $\mu_1(t)$ (orange solid line) follows the trajectory $0.3 + t$ (orange dotted line) until it is caught by $\tilde{X}_t$'s bulk (in dark gray) at $t = 0.3\,$. The same occurs for $X_t$'s spike $\lambda_1(t)$ (red solid line) and the trajectory $1 + t$ (red dotted line) until $t=1\,$.}
\label{fig:spike_spike}
\end{figure}

In this case, we consider our dynamics for $t < \min(\lambda^2, \mu^2 / q)$ so that both matrices still have a spike. In Figure \ref{fig:spike_spike}, we provide a visualisation of this situation in a specific example.

In Appendix \ref{subsec:App_spike_spike}, we demonstrate that the squared overlap $f^N(t) := \bk{\Phi_1^t}{\Psi_1^t}^2$ remains of order 1 as $N$ approaches infinity. Furthermore, it converges towards a deterministic function $f(t)$ that satisfies the differential equation
\begin{align}
\label{eq:f_equation_not_simplified}
f'(t) = \left[- \int_{\mathbb{R}} \frac{\rho(\lambda,t)}{(\lambda_1(t) - \lambda)^2}\,d\lambda - q\,\int_{\mathbb{R}} \frac{\tilde{\rho}(\mu,t)}{(\mu_1(t) - \mu)^2}\,d\mu + 2\,S\left(\lambda_1(t), \mu_1(t), t\right) \right] \, f(t) \,.
\end{align}
This expression can be simplified. Since $\lambda_1(t)$ is outside of the bulk, we have $\int_{\mathbb{R}} \frac{\rho(\lambda,t)}{\lambda_1(t) - \lambda} \, d\lambda = G\left(\lambda_1(t), t\right)\,$, which also equals $d\lambda_1 \,/\, dt\,$. Thus,
\begin{equation*}
G\left(\lambda_1(t), t\right) = \frac{1}{\lambda} \,.
\end{equation*}
Additionally, for any $z \in \mathbb{C} \setminus \mathbb{R}\,$,
\begin{align*}
- \int_{\mathbb{R}} \frac{\rho(\lambda,t)}{(z - \lambda)^2} \,d\lambda \,&= \, \partial_z G(z,t) \,=\, \frac{1}{2t}\left(1 - \frac{z}{\sqrt{z^2 - 4t}}\right) \,.
\end{align*}
Evaluating this for $z = \lambda_1(t) = \lambda + t \,/\, \lambda$ yields
\begin{equation*}
- \int_{\mathbb{R}} \frac{\rho(\lambda,t)}{(\lambda_1(t) - \lambda)^2}\,d\lambda = \frac{1}{t - \lambda^2} \,.
\end{equation*}
By performing similar calculations for $\tilde{G}\,$, we obtain
\begin{equation}
\label{eq:f_equation_simplified}
f'(t) = \left[\frac{1}{t - \lambda^2} + \frac{q}{qt - \mu^2} + \frac{2q}{\lambda \mu - qt}   \right] \,f(t)\,.
\end{equation}
This equation can be easily solved given an initial condition. At $t=0\,$, the spike of $A$ is associated with the eigenvector $\Psi_1^0 = \pm \, \Psi \,/\, \| \Psi \|$ and $\tilde{A}$'s with $\Phi_1^0 = \pm \, \Psi^{|n}\, / \, \| \Psi^{|n} \|\,$. Note that $\bk{\Psi^{|n}}{\Psi} = \| \Psi^{|n} \|^2\,$, which gives the initial condition $\bk{\Phi_1^0}{\Psi_1^0}^2 = \| \Psi^{|n} \|^2 \,/\, \| \Psi \|^2 = \mu \,/\, \lambda\,$. Finally,
\begin{equation}
\label{eq:f_solution}
f(t) = \frac{\mu}{\lambda} \, \frac{(\lambda^2 - t) \,(\mu^2 - qt)}{(\lambda \mu - qt)^2} \,.
\end{equation}
Note that this can be expressed with observable quantities only by using the relations
\begin{align*}
\lambda &= \frac{1}{2}\left(\lambda_1(t) + \sqrt{\lambda_1(t)^2 - 4t}   \right) \,, \\
\mu &= \frac{1}{2} \left(\mu_1(t) + \sqrt{\mu_1(t)^2 - 4qt}  \right) \,.
\end{align*}

\subsection{Spike-bulk overlaps}
In this second case, the truncated matrix $\tilde{X}_t$ does not have a spike, but only a semicircular bulk of radius $2 \sqrt{qt}\,$. We want to determine how its eigenvectors project onto $X_t$'s spike. Our goal is therefore to study the behaviour of $N\, \mathbb{E}\left[\bk{\Phi_i}{\Psi_1}^2\right]$ for $1 \leq i \leq n\,$, as we expect these squared overlaps to be of order $1\,/\,N\,$. To obtain an evolution equation, we introduce the following complex function:
\begin{equation*}
S_{\lambda}^{(n)}(z,t) := \sum_{i = 1}^n \, \frac{\bk{\Phi_i^t}{\Psi_1^t}^2}{z - \mu_i} \,.
\end{equation*}
It plays a role similar to the one played by $S^{(N)}$ in Section \ref{sec:Overlaps}. In Appendix \ref{subsec:App_spike_bulk}, we demonstrate that $S_{\lambda}^{(n)}$ converges towards a deterministic function $S_{\lambda}\,$, which is therefore also the limit of its mean, and satisfies the differential equation
\begin{equation}
\label{eq:S_lambda_equation}
\partial_t S_{\lambda} = -q\tilde{G}(z,t) \, \partial_z S_{\lambda} + \left(A(z,t) + B(z,t)\right)\, S_{\lambda} + C(z,t) \,,
\end{equation}
where
\begin{align*}
A(z,t) &:= \frac{2q}{\lambda (z - qt\tilde{G}(z,t)) - qt} \,, \\
B(z,t) &:= \frac{1}{t - \lambda^2} \,, \\
C(z,t) &:= \frac{q\lambda^2(z - qt\tilde{G}(z,t))}{\left(\lambda^2 - t\right)\left(\lambda(z - qt\tilde{G}(z,t)) - qt\right)^2} \,.
\end{align*}
Once more, the characteristic curve is the same as in (\ref{eq:G_tilde_implicit}). Thus, we can perform the change of variables
\begin{equation*}
\begin{cases}
y = z - qt\,\tilde{G}(z,t)  \\
s = t \,.
\end{cases}
\end{equation*}
By introducing the function $\hat{S}_{\lambda}$ defined by $S_{\lambda}(z,t) = \hat{S}_{\lambda}\left(y(z,t), s(z,t)\right)\,$, the equation simplifies to
\begin{equation*}
\partial_s \hat{S}_{\lambda} = \left(\frac{2q}{\lambda y - qs} + \frac{1}{s - \lambda^2}\right)\, \hat{S}_{\lambda} + \frac{q \lambda^2 y}{(\lambda^2 - s)(\lambda y - qs)^2} \,.
\end{equation*}
This can be solved, given that the initial condition here is $S_{\lambda}(\cdot,0) \equiv 0\,$, which implies $\hat{S}_{\lambda}(\cdot,0) \equiv 0\,$. We find
\begin{equation*}
\hat{S}_{\lambda}(y,s) = \frac{qys}{(\lambda y - qs)^2}  \,,
\end{equation*}
or, written in the original variables $z$ and $t\,$:
\begin{equation*}
S_{\lambda}(z,t) = \frac{qt\,(z - qt\,\tilde{G}(z,t))}{\left(\lambda\,\left(z - qt\,\tilde{G}(z,t)\right) - qt\right)^2} \,.
\end{equation*}
Then, if we define $g(\mu,t) := \lim_{N \to \infty} \, N\, \mathbb{E}\left[\bk{\Phi_i}{\Psi_1}^2\right]$ with $i_n \, / \, n \to x$ and $\mu = \mu(x,t)\,$,
\begin{equation*}
S_{\lambda}(z,t) = q \int_{\mathbb{R}}  \frac{g(\mu,t) \, \tilde{\rho}(\mu,t)}{z - \mu}\,d\mu \,,
\end{equation*}
so that we can use the classical Stieltjes inversion formula (\ref{eq:Stieltjes_inversion}):
\begin{equation*}
 q\,g(\mu,t)\,\tilde{\rho}(\mu,t) = \frac{1}{\pi} \, \lim_{\varepsilon \to 0^+} \, \Im\left[S_{\lambda}(\mu - i \,\varepsilon,t)\right] \,,
\end{equation*}
and obtain, after some simplifications,
\begin{equation}
\label{eq:g_solution}
g(\mu,t) = \frac{(\lambda^2 - qt)\,t}{(\lambda^2 - \lambda \mu + qt)^2} \,.
\end{equation}
The interesting thing here is that even if initially there is no information about $A$'s spike in $\tilde{A}\,$, the Brownian noise induces a coupling that leads to a non-trivial overlap between the noisy spike and $\tilde{X}_t$'s bulk. We can derive
\begin{equation*}
\lim_{N \to \infty} \, \sum_{i = 1}^n \, \mathbb{E}\left[\bk{\Phi_i^t}{\Psi_1^t}^2\right] = q \int_{\mathbb{R}} g(\mu,t)\, \tilde{\rho}(\mu,t) \, d\mu = q\,\frac{t}{\lambda^2}\,,
\end{equation*}
which shows this coupling increases linearly (just like the variance of the noise) in $t\,$, until the point where the spike is caught by the edge of the bulk at $t = \lambda^2\,$.

\begin{figure}[!h]
\centering
\includegraphics[scale=0.7]{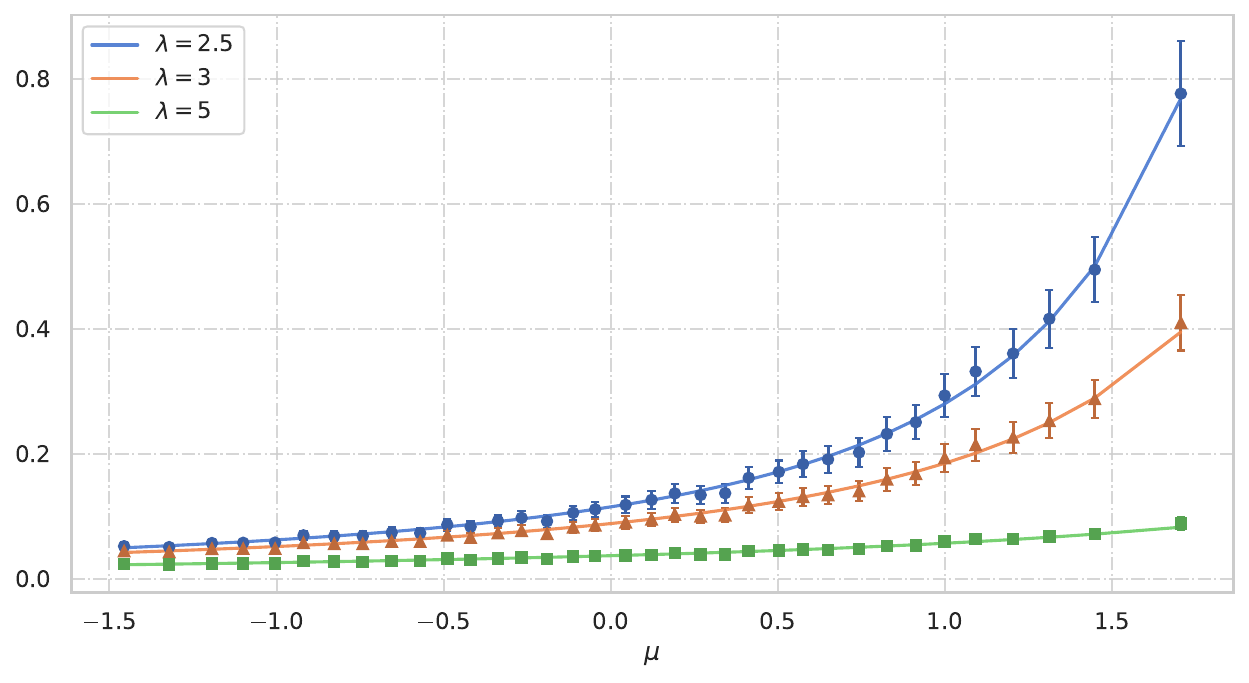}
\caption{Comparison of our theoretical formula for $g$ with numerical simulations as a function of $\mu$ in the case $N = 500\,$, $t = 1\,$, and $q = 0.7$ for different values of $\lambda = 2.5$ (blue circles for data and blue curve for theory), $\lambda = 3$ (orange triangles and orange curve) and $\lambda = 5$ (green squares and green curve). The data points are shown with 99\% confidence intervals.}
\label{fig:spike_bulk}
\end{figure}

In Figure \ref{fig:spike_bulk}, we show a comparison of this formula with simulated data of $N\,\mathbb{E}\left[\bk{\Phi_i^t}{\Psi_1^t}^2\right]$ for different values of $\lambda\,$. We used the vector $\Psi = \frac{\lambda}{\sqrt{N-n}}\,(0\,...\,0 \,1 \,...\, 1)^T$ (with $n$ zeros) to initialize the matrix $A\,$. This plot shows how, as the initial spike $\lambda$ gets closer to the edge of the spectrum (equal to $2$ in this case), the overlaps on $\tilde{X}_t$'s eigenvectors associated with large eigenvalues increase: the spike becomes more visible.

\section{Extension to the Bernoulli case}
\label{sec:Bernoulli}
The goal of this final section is to demonstrate how one can apply the results derived above in a different setup. Due to universality properties in the scaling limit in Random Matrix Theory, we expect our results to be applicable beyond the Gaussian case. Of course, this does not constitute a mathematical proof of the universality of our formulas, but rather serves as a numerical observation.

Consider $X\,$, a $N \times N$ symmetric matrix with coefficients following independent Bernoulli laws with parameter $p \in [0\,,1]\,$, rescaled by a factor $1 \,/\, \sqrt{N}\,$ , i.e. $X^{ij} = X^{ji} \sim \mathcal{B}(p) \,/\, \sqrt{N}\,$. We want to study the overlaps of its eigenvectors with those of its $n$-truncated version $\tilde{X}$ (as we defined $\tilde{X}_t$ from $X_t$ in (\ref{eq:truncation_definition})). One can view $X$ as the deterministic matrix $A$ with all coefficients equal to $p \,/\, \sqrt{N}\,$, perturbed by $H\,$, a symmetric matrix containing centered random variables with variance $p\,(1-p) \,/\, N\,$. Note that $A$ has rank 1 with its non-zero eigenvalue equal to $\lambda_N = p \,\sqrt{N}$ and associated with the eigenvector $\Psi_N = (1\,...\,1)^T \,/\, \sqrt{N}\,$. This distinguishes this application from what we have done until now: the spike eigenvalue diverges in $\sqrt{N}\,$ while the other eigenvalues form a $\mathcal{O}(1)$ bulk. Similarly, the $n$-truncated $\tilde{A}$ has rank 1 with eigenvalue $\mu_N = n \,p\, /\, \sqrt{N}$ and eigenvector $\Phi_N = (1\,...\,1\, 0\,...\,0)^T \,/\, \sqrt{n}\,$. These facts give us two intuitions:
\begin{enumerate}
\item The overlaps between eigenvectors of the bulks of $X$ and $\tilde{X}$ should asymptotically satisfy (\ref{eq:solution_noise_only}) with $t = p\,(1-p)\,$ fitting the variances of the Bernoulli variables.
\item For the spike-spike overlap, since the eigenvalues diverge in the scaling limit, one would aim to find a finite $N$ approximation of the overlap by plugging $\lambda_N$ and $\mu_N$ into (\ref{eq:f_solution}).
\end{enumerate}

\begin{figure}[h!]
	
    \centering
    \begin{minipage}[b]{0.48\textwidth}
        \centering
        \includegraphics[width = \textwidth, height = 5cm]{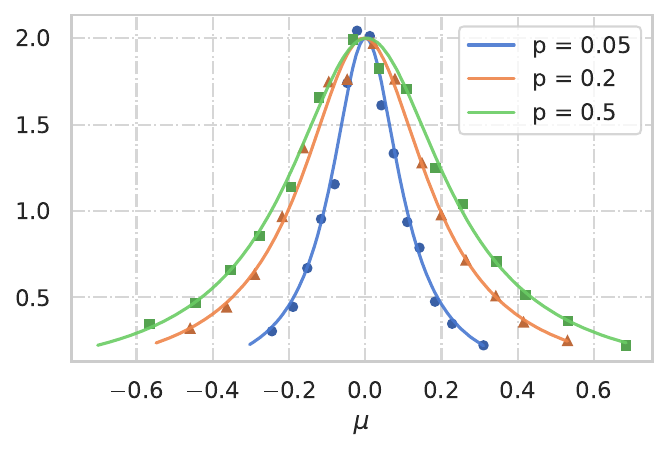}
        \label{fig:image1}
    \end{minipage}
\hspace{0cm}
    \begin{minipage}[b]{0.48\textwidth}
        \centering
        \includegraphics[width = \textwidth, height = 5cm]{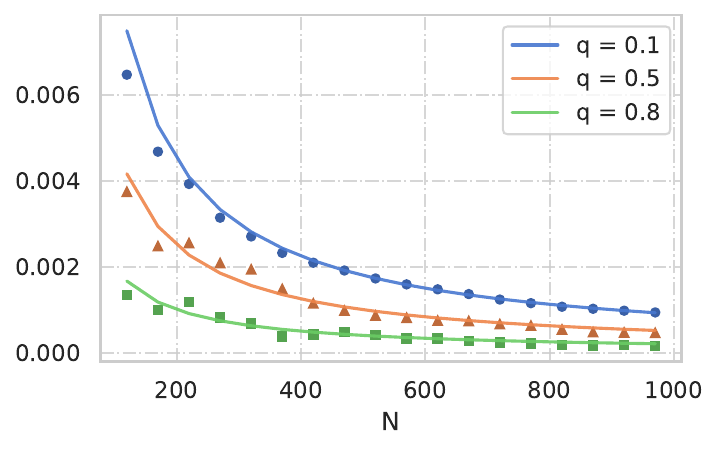}
        \label{fig:image2}
    \end{minipage}
    \caption{\textbf{Left:} comparison of formula (\ref{eq:solution_noise_only}) with simulated rescaled mean squared overlaps $N \, \mathbb{E}\left[\bk{\Phi_i}{\Psi_j}^2\right]$ of bulk eigenvectors only for $\lambda = 0$ as a function of $\mu\,$, with $N = 300$ and $q = 0.5\,$. We plot it for $p = 0.05$ (blue circles for data, blue solid curve for theory), $p = 0.2$ (orange triangles and orange solid curve) and $p = 0.5$ (green squares and green solid curve). \textbf{Right:} comparison of the $1\,/\,N$ term in (\ref{eq:bernoulli_spike_approximation}) with simulated $n\,/\,N - \mathbb{E}\left[\bk{\Phi_1^N}{\Psi_1^N}^2 \right]$ with $p = 0.7\,$, for different $N\,$. We show it for $q = 0.1$ (blue circles for data and blue solid line for theory), $q = 0.5$ (orange triangles and orange solid line) and $q = 0.7$ (green squared and green solid line).}
    \label{fig:Bernoulli}
\end{figure}

The left plot in Figure \ref{fig:Bernoulli} shows that the first point seems to hold numerically, as our theoretical formula impressively fits the bulk overlaps simulated using Bernoulli variables.

Concerning the second point, we start by noticing that in the derivation of (\ref{eq:f_solution}) (see Appendix \ref{subsec:App_spike_spike}), the largest terms we neglect are of order $1 \,/\, \sqrt{N}\,$. If we replace the order 1 eigenvalues in that context with the $\mathcal{O}(\sqrt{N})$ ones observed in the Bernoulli setup, these terms become of order $1 \,/\, N^{3/2}\,$. Therefore, this formula is valid up to this order, meaning we propose that
\begin{equation*}
\mathbb{E}\left[\bk{\Phi_1^N}{\Psi_1^N}^2 \right] = \frac{\mu_N}{\lambda_N} \, \frac{\left(\lambda_N^2 - p\,(1-p)\right) \,\left(\mu_N^2 - \frac{n}{N} \,p\, (1-p) \right)}{\left(\lambda_N \, \mu_N - \frac{n}{N} \, p\, (1-p)\right)^2} \,+\, \mathcal{O} \left(\frac{1}{N\sqrt{N}}\right) \,,
\end{equation*}
where we again took $t = p\,(1-p)$ to fit the variances of the Bernoulli variables. Then, taking the Taylor expansion after replacing the eigenvalues by their respective expressions gives
\begin{equation}
\label{eq:bernoulli_spike_approximation}
\mathbb{E}\left[\bk{\Phi_1^N}{\Psi_1^N}^2 \right] = \frac{n}{N} \,-\, \left(1 - \frac{n}{N}\right)\left(\frac{1}{p} - 1\right) \frac{1}{N} \,+\, \mathcal{O}\left(\frac{1}{N\sqrt{N}}\right) \,.
\end{equation}
This indicates that the mean squared overlap converges to $q\,$, which is also the limit of the spike-spike squared overlap of the deterministic matrices $A$ and $\tilde{A}$.
The right plot of Figure \ref{fig:Bernoulli} shows a comparison of (\ref{eq:bernoulli_spike_approximation}) up to order $1\,/\,N$  with simulated mean squared overlaps. The fit is quite remarkable. This brief computation thus demonstrates that our results can be applied to a wide range of problems, even in finite $N$ situations.

\newpage
\appendix
\section*{Appendices}
\renewcommand{\thesubsection}{\Alph{subsection}}
\renewcommand{\theequation}{\Alph{subsection}.\arabic{equation}}
\numberwithin{equation}{subsection}

\subsection{Proof of the correlation structure between the Brownian motions}
\label{subsec:correlation}
We first start by computing $c(i,l,j,k) := \bk{\Phi_i^t}{d\tilde{X}_t\,\Phi_l^t} \bk{\Psi_j^t}{dX_t \,\Psi_k^t}\,$.
Let $i\,, l \leq n$ and $j\,, k \leq N\,$:
\begin{align*}
c(i,l,j,k) &=  \left( \sum_{s, s' = 1}^N  \Phi_{is}^t \,\Phi_{ls'}^{t} \, d\tilde{X}_t^{ss'} \right) \left( \sum_{s, s' = 1}^N  \Psi_{js}^t \, \Psi_{ks'}^{t}\, dX_t^{ss'} \right) \\
&=  \left( \frac{1}{\sqrt{N}} \, \sum_{\substack{s, s' = 1 \\ s < s'}}^n \, (\Phi_{is}^t \, \Phi_{ls'}^{t} + \Phi_{is'}^{t} \, \Phi_{ls}^t) \, dB_t^{ss'} +  \sqrt{\frac{2}{N}} \,\sum_{s = 1}^n \, \Phi_{is}^t \, \Phi_{ls}^t \, dB_t^{ss} \right)  \\
& \quad \times \left( \frac{1}{\sqrt{N}} \, \sum_{\substack{s, s' = 1}}^N \,(\Psi_{js}^t \, \Psi_{ks'}^{t} + \Psi_{js'}^{t} \, \Psi_{ks}^t) \, dB_t^{ss'} +  \sqrt{\frac{2}{N}} \,\sum_{s=1}^N \, \Psi_{js}^t \, \Psi_{ks}^t \, dB_t^{ss} \right) \\
&= \frac{1}{N} \, \sum_{\substack{s, s' = 1 \\ s < s'}}^n \, (\Phi_{is}^t \, \Phi_{ls'}^{t} + \Phi_{is'}^{t}  \,\Phi_{ls}^t)  \, ( \Psi_{js}^t \, \Psi_{ks'}^{t} + \Psi_{js'}^{t} \, \Psi_{ks}^t)  \,dt \\
& \quad + \frac{2}{N}  \, \sum_{s = 1}^n\,  \Phi_{is}^t  \,\Phi_{ls}^t  \,\Psi_{js}^t \, \Psi_{ks}^t  \,dt \\
&= \frac{1}{N} \, \sum_{s, s' = 1}^n \, \left( \Phi_{is}^t\, \Phi_{ls'}^{t} \, \Psi_{js}^t \, \Psi_{ks'}^{t} + \Phi_{is}^t \, \Phi_{ls'}^{t} \, \Psi_{js'}^{t} \, \Psi_{ks}^t  \right) \, dt \,.
\end{align*}
We recall that $\Phi_{\cdot s}^t = 0$ if $s > n\,$, so that we indeed have:
\begin{equation}
\label{eq:correlation_property}
\bk{\Phi_i^t}{d\tilde{X}_t\,\Phi_l^t} \bk{\Psi_j^t}{dX_t\,\Psi_k^t} = \frac{1}{N} \left(\bk{\Phi_i^t}{\Psi_j^t} \bk{\Phi_l^t}{\Psi_k^t} + \bk{\Phi_i^t}{\Psi_k^t} \bk{\Phi_l^t}{\Psi_j^t}  \right) \, dt \,.
\end{equation}

We now aim to express the correlations not covered in Section \ref{sec:Dyson}. Projecting equation (\ref{eq:ito_correlation}) onto $\Psi_j^t\,$, we get
\begin{equation}
\label{eq:dlambda_alternative_form}
d\lambda_j^t = \bk{\Psi_j^t}{dX_t \,\Psi_j^t} + \bk{\Psi_j^t}{dX_t\, d\Psi_j^t} \,.
\end{equation}
Similarly, for $\tilde{X}_t\,$, one can derive
\begin{equation}
\label{eq:dmu_alternative_form}
d\mu_i^t = \bk{\Phi_i^t}{d\tilde{X}_t \, \Phi_i^t} + \bk{\Phi_i^t}{d\tilde{X}_t \, d\Phi_i^t} \,.
\end{equation}
This allows us to compute
\begin{align*}
dB_j(t) \, d\tilde{B}_i(t) &= \frac{N}{2} \, d\lambda_j^t \, d\mu_i^t \\
&= \frac{N}{2} \bk{\Psi_j^t}{dX_t \, \Psi_j^t} \bk{\Phi_i^t}{d\tilde{X}_t \, \Phi_i^t} \\
&= \bk{\Phi_i^t}{\Psi_j^t}^2 \, dt \,,
\end{align*}
where we used (\ref{eq:correlation_property}) in the last line. Finally, combining (\ref{eq:dW_alternative_form}) with (\ref{eq:dmu_alternative_form}) and (\ref{eq:dW_tilde_alternative_form}) with (\ref{eq:dlambda_alternative_form}):
\begin{align*}
dW_{jk}(t) \, d\tilde{B}_i(t) &= \sqrt{\frac{N}{2}} \, dW_{jk}(t) \, d\mu_i^t = \frac{N}{\sqrt{2}} \bk{\Psi_k^t}{dX_t \, \Psi_j^t} \bk{\Phi_i^t}{d\tilde{X}_t \, \Phi_i^t} \,, \\
dB_j(t) \, d\tilde{W}_{il}(t) &= \sqrt{\frac{N}{2}} \, d\lambda_j^t \, d\tilde{W}_{il}(t) = \frac{N}{\sqrt{2}} \bk{\Psi_j^t}{dX_t \, \Psi_j^t} \bk{\Phi_l^t}{d\tilde{X}_t \, \Phi_i^t} \,.
\end{align*}
Using once more the identity (\ref{eq:correlation_property}), we simplify these equalities:
\begin{equation}
\label{eq:dW_dB_correlation}
\begin{cases}
dW_{jk}(t) \, d\tilde{B}_i(t) = \sqrt{2}  \,\bk{\Phi_i^t}{\Psi_j^t} \bk{\Phi_i^t}{\Psi_k^t} \,dt \,,  \\
dB_j(t) \, d\tilde{W}_{il}(t) = \sqrt{2}\, \bk{\Phi_i^t}{\Psi_j^t} \bk{\Phi_l^t}{\Psi_j^t} \,dt \,.
\end{cases}
\end{equation}

\subsection{Method of characteristics for the Burgers equation}
\label{subsec:burgers}
We only treat the equation on $G$ since the one on $\tilde{G}$ is similar.
We introduce two functions of a new variable $s\,$ : $z(s)$ and $t(s)\,$. We also define $\hat{G}(s) := G\left(z(s),t(s)\right)\,$, so that the chain rule gives
\begin{align*}
\frac{d\hat{G}}{ds} &= \frac{dz}{ds} \, \partial_z G \left( z(s),t(s) \right) + \frac{dt}{ds} \, \partial_t G \left( z(s),t(s) \right) \\
&= \left(\frac{dz}{ds} - \hat{G} \, \frac{dt}{ds} \right) \, \partial_z G \left( z(s),t(s) \right)  \,.
\end{align*}
Therefore, if we choose $z$ and $t$ such that
\begin{equation*}
\begin{cases}
\frac{dt}{ds} = 1 \, \\
\frac{dz}{ds} = \hat{G} \,,
\end{cases}
\end{equation*}
then $d\hat{G}\,/\,ds = 0\,$. We have successfully transformed the initial equation into three simpler ones that, once solved, give
\begin{equation*}
\begin{cases}
t(s) = t(0) + s  \\
z(s) = z(0) + \hat{G}(0) \,s  \\
\hat{G}(s) = \hat{G}(0) \,.
\end{cases}
\end{equation*}
Note that $\hat{G}(0) = G\left(z(0), t(0)\right)$ and $\hat{G}(s) = G\left(z(0) + G(z,t)\,s\,, t(0) + s\right)\,$, so that when evaluating for $s = -t(0)$ and noting that $t(0)$ and $z(0)$ are free parameters, we obtain the announced implicit equation (\ref{eq:G_implicit}).

\subsection{Dynamics of the squared overlaps}
\label{subsec:App_overlaps}
Itô's formula gives
\begin{align*}
d\bk{i}{j} &= \bk{i}{dj} + \bk{di}{j} + \bk{di}{dj} \\
&= -\frac{1}{2N} \, \sum_{\substack{k = 1 \\ k \neq j}}^N  \, \frac{dt}{(\lambda_j^t - \lambda_k^t)^2} \, \bk{i}{j}  + \frac{1}{\sqrt{N}} \,\sum_{\substack{k = 1 \\ k \neq j}}^N \, \frac{dW_{jk}(t)}{\lambda_j^t - \lambda_k^t}\,\bk{i}{k} \\
&\quad  -\frac{1}{2N} \,\sum_{\substack{l = 1 \\ l \neq i}}^n\, \frac{dt}{(\mu_i^t - \mu_l^t)^2}\, \bk{i}{j}  + \frac{1}{\sqrt{N}} \,\sum_{\substack{l = 1 \\ l \neq i}}^n \, \frac{d\tilde{W}_{il}(t)}{\mu_i^t - \mu_l^t}\, \bk{l}{j} \\
&\quad + \frac{1}{N} \, \sum_{\substack{l = 1 \\ l \neq i }}^n \, \sum_{\substack{k = 1 \\ k \neq j}}^N \, \frac{\bk{i}{j}\bk{l}{k} + \bk{i}{k}\bk{l}{j}}{(\mu_i^t - \mu_l^t)(\lambda_j^t - \lambda_{k}^t)}\,\bk{l}{k}\,dt \,.
\end{align*}
Then, for the squared overlap,
\begin{align*}
d\left(\bk{i}{j}^2\right) &= 2 \,\bk{i}{j} \, d\bk{i}{j} + \left(d\bk{i}{j}\right)^2 \\
&= -\frac{1}{N} \, \sum_{\substack{k = 1 \\ k \neq j}}^N \, \frac{dt}{(\lambda_j^t - \lambda_k^t)^2}  \,\bk{i}{j}^2 + \frac{2}{\sqrt{N}}\,\sum_{\substack{k = 1 \\ k \neq j}}^N \,\frac{dW_{jk}(t)}{\lambda_j^t - \lambda_k^t}\, \bk{i}{k}\bk{i}{j} \\
&\quad -\frac{1}{N}\, \sum_{\substack{l = 1 \\ l \neq i}}^n \, \frac{dt}{(\mu_i^t - \mu_l^t)^2}\, \bk{i}{j}^2 + \frac{2}{\sqrt{N}}\, \sum_{\substack{l = 1 \\ l \neq i}}^n \, \frac{d\tilde{W}_{il}(t)}{\mu_i^t - \mu_l^t}\,\bk{l}{j}\bk{i}{j} \\
&\quad + \frac{2}{N} \, \sum_{\substack{l = 1 \\ l \neq i}}^n \, \sum_{\substack{k = 1 \\ k \neq j}}^N \, \frac{\bk{i}{j}\bk{l}{k} + \bk{i}{k}\bk{l}{j}}{(\mu_i^t - \mu_l^t)(\lambda_j^t - \lambda_{k}^t)}\, \bk{l}{k} \bk{i}{j} \,dt \\
&\quad + \frac{1}{N} \, \sum_{\substack{k = 1 \\ k \neq j}}^N \, \frac{dt}{(\lambda_j^t - \lambda_k^t)^2} \, \bk{i}{k}^2 + \frac{1}{N} \, \sum_{\substack{l = 1 \\ l \neq i}}^n \, \frac{dt}{(\mu_i^t - \mu_l^t)^2}\,\bk{l}{j}^2 \\
& \quad + \frac{2}{N} \, \sum_{\substack{l = 1 \\ l \neq i}}^n \, \sum_{\substack{k = 1 \\ k \neq j}}^N \, \frac{\bk{i}{j}\bk{l}{k} + \bk{i}{k}\bk{l}{j}}{(\mu_i^t - \mu_l^t)(\lambda_j^t - \lambda_{k}^t)} \, \bk{l}{j} \bk{i}{k} \, dt \,,
\end{align*}
which can be rearranged into
\begin{align*}
d\left(\bk{i}{j}^2\right) &= \frac{1}{N} \, \sum_{\substack{k = 1 \\ k \neq j}}^N \, \frac{\bk{i}{k}^2 - \bk{i}{j}^2}{(\lambda_j^t - \lambda_k^t)^2} \, dt + \frac{1}{N} \,  \sum_{\substack{l = 1 \\ l \neq i}}^n \, \frac{\bk{l}{j}^2 - \bk{i}{j}^2}{(\mu_i^t - \mu_l^t)^2} \,dt \\
& \quad + \frac{2}{N} \, \sum_{\substack{l = 1 \\ l \neq i}}^n \, \sum_{\substack{k = 1 \\ k \neq j}}^N \, \frac{\left( \bk{i}{j} \bk{l}{k} + \bk{i}{k} \bk{l}{j} \right)^2}{(\mu_i^t - \mu_l^t)(\lambda_j^t - \lambda_{k}^t)}\,dt \\
& \quad + \frac{2}{\sqrt{N}} \, \sum_{\substack{k = 1 \\ k \neq j}}^N \, \frac{dW_{jk}(t)}{\lambda_j^t - \lambda_k^t}\,\bk{i}{k} \bk{i}{j} + \frac{2}{\sqrt{N}} \, \sum_{\substack{l=1 \\ l \neq i}}^n \, \frac{d\tilde{W}_{il}(t)}{\mu_i^t - \mu_l^t}\,\bk{l}{j} \bk{i}{j} \,.
\end{align*}

\subsection{Limiting differential equation for the double Stieltjes transform}
\label{subsec:stieltjes_equation}

\subsubsection{Deriving the differential equation}
\label{subsubsec:deriving_S_equation}
For readability, we drop the time superscripts on the eigenvalues. By applying Itô's formula to $S^{(N)}\,$, we obtain
\begin{align*}
 dS^{(N)}(z, \tilde{z}, t) &= \frac{1}{N} \, \sum_{i=1}^n \, \sum_{j=1}^N \, \frac{d\left(\bk{i}{j}^2\right)}{(\tilde{z} - \mu_i)(z - \lambda_j)} + \frac{1}{N} \, \sum_{i=1}^n \, \sum_{j=1}^N \, \frac{\bk{i}{j}^2}{(\tilde{z} - \mu_i)^2(z-\lambda_j)}\,d\mu_i \nonumber \\
 & \quad + \frac{1}{N}  \, \sum_{i=1}^n \, \sum_{j=1}^N \,\frac{\bk{i}{j}^2}{(\tilde{z}-\mu_i)(z-\lambda_j)^2}\,d\lambda_j + \frac{1}{N} \, \sum_{i=1}^n \, \sum_{j=1}^N \,\frac{\bk{i}{j}^2}{(z - \lambda_j)(\tilde{z} - \mu_i)^3}\,(d\mu_i)^2 \nonumber \\
 &\quad + \frac{1}{N} \, \sum_{i=1}^n \, \sum_{j=1}^N \,\frac{\bk{i}{j}^2}{(z - \lambda_j)^3 (\tilde{z} - \mu_i)}\,(d\lambda_j)^2 + \frac{1}{N} \, \sum_{i=1}^n \, \sum_{j=1}^N \,\frac{\,d\left(\bk{i}{j}^2\right) \,d\mu_i}{(z-\lambda_j)(\tilde{z} - \mu_i)^2} \nonumber \\
 & \quad + \frac{1}{N} \, \sum_{i=1}^n \, \sum_{j=1}^N \,\frac{\,d\left(\bk{i}{j}^2\right)\, d\lambda_j}{(z - \lambda_j)^2(\tilde{z} - \mu_i)}  + \frac{1}{N} \, \sum_{i=1}^n \, \sum_{j=1}^N \,\frac{\bk{i}{j}^2}{(z - \lambda_j)^2 (\tilde{z} - \mu_i)^2} \,d\mu_i \, d\lambda_j \,.
\end{align*}
We need to specify certain terms:
\begin{itemize}
\item[$\bullet$] $(d\mu_i)^2 = (d\lambda_j)^2 = \frac{2}{N} \, dt\,$.
\item[$\bullet$] $d\left(\bk{i}{j}^2\right) \,d\mu_i = \frac{2\sqrt{2}}{N} \, d\tilde{B}_i \, \sum\limits_{\substack{k = 1 \\ k \neq j}}^N \, \frac{dW_{jk}}{\lambda_j - \lambda_k}\, \bk{i}{k} \bk{i}{j} = \frac{4}{N}\, \sum\limits_{\substack{k = 1 \\ k \neq j}}^N  \, \frac{\bk{i}{k}^2 \bk{i}{j}^2}{\lambda_j - \lambda_k} \,dt\,$, where we first used the independence between $\tilde{B}$ and $\tilde{W}\,$, and then (\ref{eq:dW_dB_correlation}).
\item[$\bullet$] $d\left(\bk{i}{j}^2\right) \, d\lambda_j =  \frac{2\sqrt{2}}{N} \, dB_j \, \sum\limits_{\substack{l = 1 \\ l \neq i}}^n \,\frac{d\tilde{W}_{il}}{\mu_i - \mu_l} \, \bk{l}{j} \bk{i}{j} = \frac{4}{N} \, \sum\limits_{\substack{l = 1 \\ l \neq i}}^n \, \frac{\bk{l}{j}^2 \bk{i}{j}^2}{\mu_i - \mu_l} \, dt\,$,  where we first used the independence between $B$ and $W\,$, and then (\ref{eq:dW_dB_correlation}).
\item[$\bullet$] $d\mu_i \, d\lambda_j = \frac{2}{N}  \,\bk{i}{j}^2 \, dt\,$, (see Appendix \ref{subsec:correlation}).
\end{itemize}
Thus, our Itô formula can be rewritten term by term as
\begin{align}
\label{eq:S_ito}
dS^{(N)} &= \frac{1}{N} \, \sum_{i=1}^n \, \sum_{j=1}^N \,\frac{d\left(\bk{i}{j}^2\right)}{(\tilde{z} - \mu_i)(z - \lambda_j)} + \frac{1}{N} \, \sum_{i=1}^n \, \sum_{j=1}^N \,\frac{\bk{i}{j}^2}{(\tilde{z} - \mu_i)^2(z-\lambda_j)}\,d\mu_i \nonumber \\
 & \quad + \frac{1}{N} \, \sum_{i=1}^n \, \sum_{j=1}^N \,\frac{\bk{i}{j}^2}{(\tilde{z}-\mu_i)(z-\lambda_j)^2}\,d\lambda_j + \frac{1}{N} \, \partial_{\tilde{z} \tilde{z}}^2 S^{(N)} \, dt + \frac{1}{N} \, \partial_{zz}^2 S^{(N)} \, dt \nonumber \\
 & \quad + \frac{4}{N^2} \, \sum_{i=1}^n \, \sum_{\substack{j,k=1 \\ k \neq j}}^N \,\frac{\bk{i}{k}^2 \bk{i}{j}^2}{(\lambda_j - \lambda_k)(z - \lambda_j)(\tilde{z} - \mu_i)^2}\, dt \nonumber\\
 & \quad + \frac{4}{N^2} \, \sum_{\substack{i,l=1 \\ l \neq i}}^n \, \sum_{j=1}^N \, \frac{\bk{l}{j}^2 \bk{i}{j}^2}{(\mu_i - \mu_l)(z-\lambda_j)^2(\tilde{z} - \mu_i)}\, dt + \frac{2}{N^2} \, \sum_{i=1}^n \, \sum_{j=1}^N \, \frac{\bk{i}{j}^4}{(z - \lambda_j)^2(\tilde{z} - \mu_i)^2} \, dt \,. 
\end{align}
Notice that, since the squared overlaps are expected to be of order $1\,/\,N\,$, all terms vanish except the first three, which remain of order 1 in the scaling limit. Therefore, we are interested in combining these remaining terms to identify functions of $S^{(N)}$ and of its derivatives. First, we decompose these terms and introduce some new notations. Let us define:
\begin{itemize}
\item[$\bullet$] $d\Sigma_{\mu} := \frac{1}{N} \, \sum\limits_{i=1}^n \, \sum\limits_{j=1}^N \, \frac{\bk{i}{j}^2}{(\tilde{z} - \mu_i)^2 (z- \lambda_j)}\,d\mu_i\,$,
\item[$\bullet$] $d\Sigma_{\lambda} := \frac{1}{N} \, \sum\limits_{i=1}^n \, \sum\limits_{j=1}^N \, \frac{\bk{i}{j}^2}{(\tilde{z} - \mu_i)(z - \lambda_j)^2} \,d\lambda_j\,$.
\end{itemize}
The first term in (\ref{eq:S_ito}), involving $d\left(\bk{i}{j}^2\right)\,$, needs to be decomposed using our previous Itô formula on the squared overlaps (\ref{eq:ito_squared_overlaps}). It is equal to
\begin{equation*}
\left(I_{\lambda} + I_{\mu} + I_{\lambda \mu}\right) \, dt + dI_W + dI_{\tilde{W}} \,,
\end{equation*}
where:
\begin{itemize}
\item[$\bullet$] $I_{\lambda} := \frac{1}{N^2} \, \sum\limits_{i=1}^n \, \sum\limits_{\substack{j,k = 1 \\ k \neq j}}^N \, \frac{\bk{i}{k}^2 - \bk{i}{j}^2}{(\lambda_j - \lambda_k)^2(z - \lambda_j)(\tilde{z} - \mu_i)}$ contains interactions between the $\lambda_j\,$.
\item[$\bullet$] $I_{\mu} := \frac{1}{N^2} \, \sum\limits_{\substack{i,l = 1\\ l \neq i}}^n \, \sum\limits_{\substack{j = 1}}^N \, \frac{\bk{l}{j}^2 - \bk{i}{j}^2}{(\mu_i - \mu_l)^2 (\tilde{z} - \mu_i)(z - \lambda_j)}$ contains interactions between the $\mu_i\,$.
\item[$\bullet$] $I_{\lambda \mu} := \frac{2}{N^2} \, \sum\limits_{\substack{i,l = 1 \\ l \neq i}}^n \sum\limits_{\substack{j,k = 1 \\ k \neq j}}^N \, \frac{\left( \bk{i}{j} \, \bk{l}{k} + \bk{i}{k} \, \bk{l}{j} \right)^2}{(\mu_i - \mu_l)(\lambda_j - \lambda_{k})(\tilde{z} - \mu_i)(z - \lambda_j)}$ contains both interactions between the $\lambda_j$ and between the $\mu_i\,$.
\end{itemize}
We also denoted by $dI_W$ and $dI_{\tilde{W}}$ the terms involving the Brownian motions:
\begin{itemize}
\item[$\bullet$] $dI_W := \frac{2}{N \sqrt{N}} \, \sum\limits_{i=1}^n \sum\limits_{\substack{j,k = 1 \\ k \neq j}}^N \, \frac{dW_{jk} \, \bk{i}{k} \bk{i}{j}}{(\lambda_j - \lambda_k)(z - \lambda_j)(\tilde{z} - \mu_i)} \,$,
\item[$\bullet$] $dI_{\tilde{W}} := \frac{2}{N \sqrt{N}} \, \sum\limits_{\substack{i,l = 1 \\ l \neq i}}^n \sum\limits_{j=1}^N \frac{d\tilde{W}_{il}  \, \bk{l}{j} \bk{i}{j}}{(\mu_i - \mu_l)(\tilde{z}-\mu_i)(z - \lambda_j)}\,$.
\end{itemize}
Let us now detail the steps of our computation:
\begin{enumerate}
\item First, we demonstrate that the terms involving the Brownian motions almost surely converge to zero as $N$ tends to infinity. This is a key step, as it explicitly shows how the randomness in $S^{(N)}$ vanishes in the scaling limit, leaving only terms that can be related to itself or to its derivatives.
\item Next, we establish that by combining $\Sigma_{\lambda}$ and $I_{\lambda}\,$, we can identify terms involving $\partial_z S^{(N)}\,$. Similarly, the combination of $\Sigma_{\mu}$ with $I_{\mu}$ will yield the $\partial_{\tilde{z}}$ part of the equation.
\item Finally, we manipulate the term $I_{\lambda \mu}$ to show that it almost surely converges to $S^2\,$.
\end{enumerate}
To manipulate the sums present in our expressions, we will use the following identity:
\begin{equation}
\tag{I} 
\label{eq:identity}
\frac{1}{(a - b)\,(z - a)} - \frac{1}{(a-b)\,(z-b)} = \frac{1}{(z-a)\,(z-b)} \,.
\end{equation}
We will also make use of two types of symmetrisations:
\begin{itemize}
\item[$\bullet$] \underline{Symmetrisation 1:} If $a_{kl} = a_{lk}\,$,
\begin{equation}
\tag{S1} 
\label{eq:S1}
\sum_{\substack{k,l \\ k \neq l}} \, \frac{a_{kl}}{(b_k - b_l) \, (z - b_k)} = \frac{1}{2} \, \sum_{\substack{k,l \\ k \neq l}} \, \frac{a_{kl}}{(z - b_k)\, (z - b_l)} \,.
\end{equation}
To prove this, we copy the sum $S$ into $S/2 + S/2\,$, invert the indices of the second sum, and apply the identity (\ref{eq:identity}).

\item[$\bullet$] \underline{Symmetrisation 2:}
\begin{equation}
\tag{S2} 
\label{eq:S2}
\sum_{\substack{k,l \\ k \neq l}}\, a_{kl} + a_{lk} = 2\,  \sum_{\substack{k,l \\ k \neq l}} \, a_{kl} \,.
\end{equation}
This relation is derived by expanding the sum, inverting the indices in the second term, and combining the two sums.
\end{itemize}
In this subsection, we will explicitly indicate when one of these properties is used and the concerned indices.

We recall that we are working with a fixed interval $[t\,, t + dt]\,$, and with fixed $z\,$, $\tilde{z} \in \mathbb{C} \setminus \mathbb{R}\,$, all of which are independent of $N\,$.

\paragraph{Brownian terms.}

We show here that $dI_W \overset{\text{a.s.}}{\rightarrow} 0$ as $N$ goes to infinity. The method is identical for $dI_{\tilde{W}}\,$. The sum $dI_W$ has mean zero due to the independence between $dW(t)$ and $\Phi^t$, $\Psi^t$, $\mu^t$, $\lambda^t$. We are going to use some manipulations to demonstrate that its variance is of order $\mathcal{O}(1 \,/\, N^2)\,$. Therefore, by applying Borel-Cantelli's lemma to the events $\{ |dI_W| > \varepsilon \}\,$, with their probabilities bounded from above by Bienaymé-Tchebychev's inequality, we obtain almost sure convergence.

First, we apply symmetrisation (\ref{eq:S1}) to the indices $j$ and $k\,$, with $a_{jk} = dW_{jk}\,\bk{i}{k}\,\bk{i}{j}$ and $b_j = \lambda_j\,$. This transforms our random variable into
\begin{equation*}
dI_W = \frac{1}{N\sqrt{N}} \, \sum_{i=1}^n \sum_{\substack{j,k = 1 \\ k \neq j}}^N \, \frac{dW_{jk} \, \bk{i}{k} \bk{i}{j}}{(z - \lambda_j)(z- \lambda_k)(\tilde{z}- \mu_i)} \,.
\end{equation*}
This prevents having a diverging term $\frac{1}{N} \, \sum_{j \neq k} \, \frac{1}{(\lambda_j - \lambda_k)^2}$ in the variance, as all the eigenvalues lie in the bulk. Note that if we only account for the fact that the overlaps are of order $1 \,/\, \sqrt{N}\,$, this sum would be of order $\mathcal{O}(\sqrt{N})\,$, but since we are summing many independent variables (the $dW$), there is a concentration of measure, and the variance will initially appear to be of order $1\,/\,N\,$. The variance reads:
\begin{align*}
\mathbb{E}\left[ \left|dI_W \right|^2 \right] &= \frac{1}{N^3} \, \sum_{\substack{i,l=1}}^n \, \sum_{\substack{j,k,j',k'=1 \\ k \neq j \\ k' \neq j'}}^N \, \mathbb{E}\left[ \frac{ dW_{jk}\,dW_{j'k'} \, \bk{i}{k}  \bk{i}{j} \bk{l}{k'} \bk{l}{j'} }{(z - \lambda_j)(z^* - \lambda_{j'})(z - \lambda_k)(z^* - \lambda_{k'})(\tilde{z} - \mu_i)(\tilde{z}^* - \mu_{l})} \right] \,,
\end{align*}
where $z^*$ stands for the complex conjugate of $z\,$. Note that $dW_{jk}\,  dW_{j'k'} \neq 0$ only if $k' = k$ and $j' = j$ or if $k' = j$ and $j' = k\,$. In both cases $\bk{l}{k'} \bk{l}{j'} = \bk{l}{k} \bk{l}{j}\,$, thus,
\begin{align*}
\mathbb{E}\left[ \left|dI_W \right|^2 \right] &= \frac{2}{N^3} \,  \sum_{\substack{i,l=1}}^n \, \sum_{\substack{j,k=1 \\ k \neq j}}^N \, \mathbb{E}\left[ \frac{\bk{i}{k} \bk{i}{j} \bk{l}{k} \bk{l}{j} }{|z - \lambda_j|^2|z - \lambda_k|^2 (\tilde{z} - \mu_i)(\tilde{z}^* - \mu_{l})} \right] \, dt \,.
\end{align*}
This sum is indeed of order $1\,/\,N\,$. By exploiting its specific structure, we can increase this bound to $1\,/\,N^2\,$, ensuring almost sure convergence. Similar manipulations will be applied several times in the remainder of this paper, hence we present the method here as a general property.
The key observation here is that the variance can be rewritten as
\begin{equation}
\label{eq:variance_brownian_term}
\mathbb{E}\left[ \left|dI_W \right|^2 \right] =  \frac{2}{N^3} \, \sum_{\substack{j,k=1 \\ k \neq j }}^N \, \mathbb{E}\left[ \frac{1}{|z - \lambda_j|^2|z-\lambda_k|^2} \,  \left| \sum_{\substack{i=1}}^n \, \frac{\bk{i}{k} \bk{i}{j}}{\tilde{z} - \mu_i} \right|^2 \right] \, dt \,.
\end{equation}
We now introduce the following property:

\begin{itemize}
\item[$\bullet$] \underline{Reduction Property:}

For $z\,, \tilde{z} \in \mathbb{C} \setminus \mathbb{R}$ and for any $p \geq 0\,$, we have
\begin{equation}
\tag{R}
\label{eq:property_R_bulk}
\frac{1}{N^2} \,  \sum_{\substack{j,k = 1 \\ k \neq j}}^N \, \frac{1}{|z - \lambda_j|^p |z - \lambda_k|^p} \, \left| \sum_{i=1}^n \frac{\bk{i}{j}\bk{i}{k}}{\tilde{z} - \mu_i} \right|^2 = \mathcal{O}\left(\frac{1}{N}\right) \,.
\end{equation}

\noindent
\textit{Proof}

Let us denote the sum under investigation by $\Sigma\,$. Using the upper bound $1 \,/\, |z - \lambda|^2 \leq 1 \,/\, \Im(z)^2\,$, we have
\begin{align*}
\Sigma &\leq \frac{1}{N^2 \,\Im(z)^{2p}} \,  \sum_{\substack{j,k = 1 \\ k \neq j}}^N \, \left| \sum_{\substack{i=1}}^n \, \frac{\bk{i}{j} \bk{i}{k}}{\tilde{z} - \mu_i} \right|^2  \\
&\leq \frac{1}{N^2 \,\Im(z)^{2p}} \,  \sum_{\substack{j,k = 1 \\ k \neq j}}^N  \, \sum_{\substack{i,l = 1}}^n \, \frac{\bk{i}{j} \bk{i}{k} \bk{l}{j} \bk{l}{k}}{(\tilde{z} - \mu_i)(\tilde{z}^* - \mu_{l})} \\
& \leq \frac{1}{N^2 \,\Im(z)^{2p}} \, \sum_{\substack{j=1}}^N  \, \sum_{\substack{i,l = 1}}^n \, \frac{\bk{i}{j} \bk{l}{j}}{(\tilde{z} - \mu_i)(\tilde{z}^* - \mu_{l})} \, \sum_{\substack{k=1\\ k \neq j}}^N \,  \bk{i}{k} \bk{l}{k} \,.
\end{align*}
Since $\Psi_1^t\,,...\,, \Psi_N^t$ form an orthonormal basis of $\mathbb{R}^N$, we recognize that
\begin{equation*}
\sum_{\substack{k = 1 \\ k \neq j}}^N \, \bk{i}{k} \bk{l}{k} = \bk{\Phi_i^t}{\Phi_{l}^t} - \bk{i}{j} \bk{l}{j} = \delta_{il} - \bk{i}{j} \bk{l}{j} \,.
\end{equation*}
Therefore, our inequality becomes
\begin{align*}
\Sigma &\leq  \frac{1}{N^2 \,\Im(z)^{2p}} \,  \left( \sum_{i=1}^n \, \sum_{j=1}^N \, \frac{\bk{i}{j}^2}{|\tilde{z} - \mu_i|^2} - \sum_{\substack{i,l = 1}}^n \, \sum_{j = 1}^N \, \frac{\bk{i}{j}^2 \bk{l}{j}^2}{(\tilde{z} - \mu_i)(\tilde{z}^* - \mu_{l})}   \right) \\
&\leq \frac{1}{N^2 \,\Im(z)^{2p}} \, \left( \sum_{\substack{i = 1}}^{n} \, \frac{1}{|\tilde{z} - \mu_i|^2} - \sum_{j=1}^N \, \left| \sum_{\substack{i = 1}}^n \, \frac{\bk{i}{j}^2}{\tilde{z} - \mu_i}  \right|^2  \right) \\
&\leq \frac{1}{\Im(z)^{2p} \,\Im(\tilde{z})^2} \,  \frac{n}{N^2} \,,
\end{align*}
by applying the same complex upper bound to the first sum.
Since $n\,/\,N \to q\,$, we indeed have $\Sigma = \mathcal{O}(1 \,/\, N)\,$, which concludes our proof.
\end{itemize}

\noindent
Applying (\ref{eq:property_R_bulk}) to (\ref{eq:variance_brownian_term}) yields that the variance is of order $1\,/\,N^2\,$.

\paragraph{Partial derivative terms.}

Our goal here is to show the following convergence:
\begin{equation}
\label{eq:limit_lambda_part}
I_{\lambda}\, dt + d\Sigma_{\lambda} \longrightarrow -G(z,t) \,  \partial_z S(z, \tilde{z}, t)\,dt \,.
\end{equation}
We demonstrate our method for the $\lambda$-case, and similar manipulations lead to
\begin{equation*}
I_{\mu}\,dt + d\Sigma_{\mu} \longrightarrow -q \tilde{G}(\tilde{z},t) \, \partial_{\tilde{z}} S(z, \tilde{z}, t) \, dt \,.
\end{equation*}
We begin by manipulating $I_{\lambda}\,$. By applying symmetrisation (\ref{eq:S1}) to the indices $j$ and $k\,$, with $a_{jk} = \frac{\bk{i}{k}^2 - \bk{i}{j}^2}{(\lambda_j - \lambda_k)}$ and $b_j = \lambda_j\,$, we obtain
\begin{equation*}
I_{\lambda} = \frac{1}{2N^2} \, \sum_{i=1}^n \, \sum_{\substack{j,k = 1 \\ k \neq j}}^N \, \frac{\bk{i}{k}^2 - \bk{i}{j}^2}{(\lambda_j - \lambda_k)(z - \lambda_j)(z- \lambda_k)(\tilde{z} - \mu_i)}  \,,
\end{equation*}
which we can split and regroup as in (\ref{eq:S2}):
\begin{equation*}
I_{\lambda} = \frac{1}{N^2} \, \sum_{i=1}^n \, \sum_{\substack{j,k = 1 \\ k \neq j}}^N \, \frac{\bk{i}{k}^2}{(\lambda_j - \lambda_k)(z - \lambda_j)(z- \lambda_k)(\tilde{z} - \mu_i)} \,.
\end{equation*}
Finally, applying identity (\ref{eq:identity}) with $a = \lambda_j$ and $b = \lambda_k\,$, we get
\begin{align*}
I_{\lambda} &= \overbrace{\frac{1}{N^2}\, \sum_{i=1}^n \, \sum_{\substack{j,k = 1 \\ k \neq j}}^N \, \frac{\bk{i}{k}^2}{(\lambda_j - \lambda_k)(z - \lambda_k)^2(\tilde{z} - \mu_i)}}^{A_{\lambda}} + \frac{1}{N^2} \, \sum_{i=1}^n \, \sum_{\substack{j,k = 1 \\ k \neq j}}^N \, \frac{\bk{i}{k}^2}{(z - \lambda_k)^2(z - \lambda_j)(\tilde{z} - \mu_i)} \\
&= A_{\lambda} + \left(\frac{1}{N} \,  \sum_{j = 1}^N \,  \frac{1}{z - \lambda_j} \right) \left( \frac{1}{N} \,  \sum_{i=1}^n \, \sum_{k = 1}^N \, \frac{\bk{i}{k}^2}{(z - \lambda_k)^2(\tilde{z} - \mu_i)} - \frac{1}{N} \, \frac{\bk{i}{j}^2}{(z - \lambda_j)^2(\tilde{z} - \mu_i)} \right) \\
&= A_{\lambda} - G_N(z,t) \,  \partial_z S^{(N)} - \frac{1}{2N} \, \partial_{zz}^2 S^{(N)} \,.
\end{align*}
Next, we observe that $A_{\lambda}$ can be cancelled by $d\Sigma_{\lambda}\,$. Indeed,
\begin{align*}
d\Sigma_{\lambda} &= \frac{\sqrt{2}}{N \sqrt{N}} \,  \sum_{i=1}^n \, \sum_{j=1}^N \, \frac{\bk{i}{j}^2}{(z - \lambda_j)^2(\tilde{z} - \mu_i)}\,dB_j + \frac{1}{N^2} \, \sum_{i=1}^n \, \sum_{\substack{j,k=1 \\ k \neq j}}^N \, \frac{\bk{i}{j}^2}{(\lambda_j - \lambda_k)(z- \lambda_j)^2(\tilde{z} - \mu_i)} \, dt \\
&\hspace{-1.1cm} \stackrel{\text{(Inverting $j$ and $k$)}}{=} \frac{\sqrt{2}}{N \sqrt{N}} \, \sum_{i=1}^n \,  \sum_{j=1}^N \, \frac{\bk{i}{j}^2}{(z - \lambda_j)^2(\tilde{z} - \mu_i)}\,dB_j - A_{\lambda} \, dt \,,
\end{align*}
so that
\begin{equation*}
I_{\lambda}\, dt + \Sigma_{\lambda} = -G_N(z,t) \,  \partial_z S^{(N)} \, dt - \frac{1}{2N} \,  \partial_{zz}^2 S^{(N)} \, dt + \frac{\sqrt{2}}{N \sqrt{N}} \, \sum_{i=1}^n \, \sum_{j=1}^N \, \frac{\bk{i}{j}^2}{(z - \lambda_j)^2(\tilde{z} - \mu_i)}\,dB_j \,,
\end{equation*}
where only the first term remains of order 1 in the scaling limit, while all other go to $0\,$. This gives the announced convergence (\ref{eq:limit_lambda_part}).

\paragraph{Square term.}

Let us first expand the square in $I_{\lambda \mu}\,$:
\begin{align*}
I_{\lambda \mu} &= \frac{2}{N^2} \,  \sum_{\substack{i,l=1 \\ l \neq i}}^n \, \sum_{\substack{j,k = 1 \\ k \neq j}}^N \, \frac{\bk{i}{j}^2 \bk{l}{k}^2 + \bk{i}{k}^2 \bk{l}{j}^2}{(\mu_i - \mu_l)(\lambda_j - \lambda_{k})(\tilde{z}- \mu_i)(z - \lambda_j)} \\
& \quad + \frac{4}{N^2} \,  \sum_{\substack{i,l=1 \\ l \neq i}}^n \,  \sum_{\substack{j,k=1\\ k \neq j}}^N \, \frac{\bk{i}{j} \bk{l}{k} \bk{i}{k} \bk{l}{j}}{(\mu_i - \mu_l)(\lambda_j - \lambda_{k})(\tilde{z}- \mu_i)(z - \lambda_j)} \,.
\end{align*}
We begin with the first sum. By applying symmetrisation (\ref{eq:S1}) to $(i,l)$ and then to $(j,k)\,$, it becomes
\begin{equation*}
\frac{1}{2N^2} \, \sum_{\substack{i,l=1 \\ l \neq i}}^n \, \sum_{\substack{j,k = 1 \\ k \neq j}}^N \,  \frac{\bk{i}{j}^2 \bk{l}{k}^2 + \bk{i}{k}^2 \bk{l}{j}^2}{(\tilde{z} - \mu_i)(\tilde{z}- \mu_l)(z - \lambda_j)(z - \lambda_{k})} \,.
\end{equation*}
Then, symmetrisation (\ref{eq:S2}) transforms it into
\begin{equation*}
\frac{1}{N^2} \,  \sum_{\substack{i,l=1 \\ l \neq i}}^n \,  \sum_{\substack{j,k = 1 \\ k \neq j}}^N \, \frac{\bk{i}{j}^2 \bk{l}{k}^2}{(\tilde{z} - \mu_i)(\tilde{z}- \mu_l)(z - \lambda_j)(z - \lambda_{k})} \,.
\end{equation*}
We observe that if we add the diagonal terms ($l=i$ and $k = j$), we exactly obtain $\left(S^{(N)}\right)^2$. One can verify that the latter all go to 0 as $N$ goes to infinity, so that this sum converges to $S^2$.

Finally, we address the second sum in $I_{\lambda \mu}$ and show it converges to 0 as $N \to \infty\,$. We first symmetrize it twice, applying (\ref{eq:S1}) to $(i,l)$ and to $(j,k)\,$, and obtain
\begin{equation*}
\frac{1}{N^2}  \, \sum_{\substack{i,l=1 \\ l \neq i}}^n \, \sum_{\substack{j,k = 1 \\ k \neq j}}^N \frac{\bk{i}{j} \bk{l}{k} \bk{i}{k}\bk{l}{j}}{(\tilde{z} - \mu_i)(\tilde{z} - \mu_l)(z - \lambda_j)(z - \lambda_{k})} \,,
\end{equation*}
which can be rewritten as (adding the diagonal terms $l = i$ that go to $0$ in the scaling limit)
\begin{equation*}
\frac{1}{N^2} \, \sum_{\substack{j,k=1 \\ k \neq j}}^N \, \frac{1}{(z - \lambda_j)(z - \lambda_{k})} \, \left(\sum_{\substack{i=1}}^n \, \frac{\bk{i}{j}\bk{i}{k}}{\tilde{z} - \mu_i}  \right)^2 \,.
\end{equation*}
Thus, its modulus is bounded from above by
\begin{equation*}
\frac{1}{N^2} \,  \sum_{\substack{j,k = 1 \\ k \neq j}}^N \, \frac{1}{|z - \lambda_j||z - \lambda_{k}|} \,  \left| \sum_{i=1}^n \, \frac{\bk{i}{j}\bk{i}{k}}{\tilde{z} - \mu_i} \right|^2 \,,
\end{equation*}
which is of order $1\,/\,N$ using the reduction property (\ref{eq:property_R_bulk}). Therefore, the sum goes to 0 as $N \to \infty\,$.

We can now conclude by regrouping all our terms. We have proven that the limit function $S$ almost surely satisfies the equation (\ref{eq:S_equation}):
\begin{equation*}
\partial_t S = -G(z,t) \, \partial_z S -q\tilde{G}(\tilde{z},t) \, \partial_{\tilde{z}} S + S^2 \,.
\end{equation*}

\subsubsection{Solving the differential equation}

We solve this equation using the method of characteristics. We introduce three functions of a new variable $s\,$: $z(s)\,, \tilde{z}(s)\,, t(s)\,$. We also define $\hat{S}(s) := S\left(z(s), \tilde{z}(s), t(s)\right)\,$. The chain rule gives
\begin{align*}
\frac{d\hat{S}}{ds} &= \partial_z S \, \frac{dz}{ds} + \partial_{\tilde{z}}S \,  \frac{d\tilde{z}}{ds} + \partial_tS \, \frac{dt}{ds} \\
&= \partial_z S \, \left( \frac{dz}{ds} - G\left(z(s), t(s)\right) \, \frac{dt}{ds}  \right) + \partial_{\tilde{z}} S \, \left( \frac{d\tilde{z}}{ds} - q\tilde{G}\left(\tilde{z}(s), t(s)\right)\, \frac{dt}{ds} \right) + \hat{S}^2 \, \frac{dt}{ds} \,.
\end{align*}
Therefore, if we choose our three functions such that
\begin{equation*}
\begin{cases}
\frac{dz}{ds} = G\left(z(s),t(s)\right)  \\
\frac{d\tilde{z}}{ds} = q \tilde{G}\left(\tilde{z}(s), t(s)\right)  \\
\frac{dt}{ds} = 1 \,,
\end{cases}
\end{equation*}
then we will have $\frac{d\hat{S}}{ds} = \hat{S}^2\,$, i.e.
\begin{equation}
\hat{S}(s) = \frac{\hat{S}(0)}{1 - s \, \hat{S}(0)} \,.
\label{eq:ricatti_sol}
\end{equation}
Also, we know that $G$ satisfies $\partial_t G = -G \, \partial_z G$ and $\tilde{G}$ satisfies $ \partial_t \tilde{G} = - q \tilde{G}\, \partial_z \tilde{G}\,$, which gives us $\frac{d^2z}{ds^2} = \frac{d^2\tilde{z}}{ds^2} = 0\,$. Therefore, we have
\begin{equation*}
\begin{cases}
z(s) = G\left(z(0), t(0)\right) s + z(0) \\
\tilde{z}(s) = q \tilde{G}\left(\tilde{z}(0), t(0)\right) s + \tilde{z}(0)  \\
t(s) = s + t(0) \,.
\end{cases}
\end{equation*}
Substituting this into (\ref{eq:ricatti_sol}), evaluating for $s = -t(0)\,$, and noticing that $z(0)\,,\, \tilde{z}(0)\,\,,\, t(0)$ are free parameters that we can choose, we obtain the announced explicit solution (\ref{eq:S_general_solution}):
\begin{equation*}
S(z, \tilde{z}, t) = \frac{S\left(z - t\,G(z,t)\,,\, \tilde{z} - qt\,\tilde{G}(\tilde{z}, t)\,,\, 0\right)}{1 - t\, S\left(z - t\,G(z,t)\,,\, \tilde{z} - qt\,\tilde{G}(\tilde{z}, t)\,,\, 0\right)} \,.
\end{equation*}

\subsection{Closed-form formula in the case of a pure noise matrix}
\label{subsec:cauchy_simplification}
In the case $A \equiv 0\,$, we have (\ref{eq:S_solution_noise_only}):
\begin{equation*}
S(z, \tilde{z}, t) = \frac{q}{(z - tG(z,t)) \, \left(\tilde{z} - qt \tilde{G}\left(\tilde{z},t\right)\right) - qt} \,.
\end{equation*}
For fixed $\lambda \in [-2 \sqrt{t}\,, 2 \sqrt{t}]$ and $\mu \in [-2 \sqrt{qt}\,, 2 \sqrt{qt}]\,$, we introduce
\begin{align*}
S_{\pm} &:= \lim_{\varepsilon \to 0^+} S(\lambda - i\,\varepsilon, \mu \pm i \,\varepsilon, t) \\
&= \frac{q}{\left(\lambda - t\,v(\lambda,t) - i \,\pi t \,\rho(\lambda,t)\right) \,\left(\mu - qt \,\tilde{v}(\mu,t) \pm i \,q \pi t\, \tilde{\rho}(\mu,t)\right) - qt} \,,
\end{align*}
so that by using (\ref{eq:double_Stieltjes_inversion}),
\begin{equation*}
W(\mu, \lambda , t) = \frac{1}{2 \pi^2 q \rho(\lambda,t) \tilde{\rho}(\mu,t)} \, \Re(S_+ - S_-) \,.
\end{equation*}
We use the following notations:
\begin{itemize}
\item[$\bullet$] $V := \lambda - t\,v(\lambda,t)\,$. Knowing that $v(\lambda,t) = \lambda \,/\, 2t\,$, we have $ V = \lambda \,/\, 2\,$.
\item[$\bullet$] $R := \pi t\, \rho(\lambda,t)\,$.
\item[$\bullet$] $\tilde{V} := \mu - qt \,\tilde{v}(\mu,t) = \mu \,/\, 2\,$.
\item[$\bullet$] $\tilde{R} := q \pi t \,\tilde{\rho}(\mu,t)\,$.
\end{itemize}
We begin by simplifying
\begin{align*}
S_{\pm} &= \frac{q}{(V - iR)(\tilde{V} \pm i\tilde{R}) -qt}\\
&= \frac{q}{V \tilde{V} \pm R \tilde{R} - qt -i (R\tilde{V} \mp \tilde{R}V)} \\
&= q \frac{V \tilde{V} \pm R \tilde{R} - qt + i(R\tilde{V} \mp \tilde{R}V)}{(V \tilde{V} \pm R \tilde{R} - qt)^2 + (R\tilde{V} \mp \tilde{R}V)^2} \,.
\end{align*}
This allows us to identify its real part. Our goal is to obtain a simplified form of
\begin{equation*}
\Re(S_+ - S_-) = q\frac{V \tilde{V} + R \tilde{R} - qt}{(V \tilde{V} + R \tilde{R} - qt)^2 + (R\tilde{V} - \tilde{R}V)^2} - q\frac{V \tilde{V} - R \tilde{R} - qt}{(V \tilde{V} - R \tilde{R} - qt)^2 + (R\tilde{V} + \tilde{R}V)^2} \,.
\end{equation*}
We denote it as $\Re(S_+ - S_-) = N_+ \,/\, D_+ \,-\, N_- \,/\, D_-\,$. We combine the two fractions under the product of the denominators and start by factorizing the numerator:
\begin{align*}
N_+ \,D_- - N_- \,D_+ &=  q\left(V \tilde{V} + R \tilde{R} - qt\right)\left((V \tilde{V} - R \tilde{R} - qt)^2 + (R\tilde{V} + \tilde{R}V)^2\right)\\
& \quad  - q\left(V \tilde{V} - R \tilde{R} - qt \right) \left((V \tilde{V} + R \tilde{R} - qt)^2 + (R\tilde{V} - \tilde{R}V)^2 \right) \,.
\end{align*}
Using the identities $(x + y)^2 - (x-y)^2 = 4xy$ and $(x+y)^2 + (x-y)^2 = 2x^2 + 2y^2\,$, we transform it into
\begin{align*}
N_+\,D_- - N_-\,D_+ &= 4q\left(V\tilde{V} -qt\right)\left(-(V\tilde{V} -qt)R\tilde{R} + V\tilde{V}R\tilde{R} \right) \\
& \quad + 2qR\tilde{R} \left( \left(V\tilde{V}-qt\right)^2 + R^2 \tilde{R}^2 + R^2 \tilde{V}^2 + \tilde{R}^2 V^2  \right) \\
&= 2qR\tilde{R} \left(V^2 \tilde{V}^2 + R^2 \tilde{R}^2 + R^2 \tilde{V}^2 + \tilde{R}^2 V^2 -q^2 t^2   \right) \\
&= 2qR\tilde{R} \left[ \left(V^2 + R^2 \right)\left(\tilde{V}^2 + \tilde{R}^2 \right) -q^2t^2 \right] \,.
\end{align*}
Now we recall that
\begin{equation*}
R^2 = \pi^2 t^2 \, \frac{4t - \lambda^2}{4 \pi^2 t^2} = t - \frac{\lambda^2}{4} \,,
\end{equation*}
so that $V^2 + R^2 = t\,$. Similarly, $\tilde{V}^2 + \tilde{R}^2 = qt\,$. Finally, our numerator equals $2q^2(1-q)t^2R\tilde{R}\,$.

We finish by simplifying the denominator, which is the product of the two initial denominators. First, the denominators of each fraction are
\begin{align*}
D_{\pm} &= (V\tilde{V} - qt \pm R\tilde{R})^2 + (R\tilde{V} \mp \tilde{R}V)^2 \\
&= (V^2 + R^2)(\tilde{V}^2 + \tilde{R}^2) + q^2t^2 -2qtV\tilde{V} \mp 2qtR\tilde{R} \\
&= qt\left((1+q)t - 2V\tilde{V} \mp 2R\tilde{R} \right) \,.
\end{align*}
Thus, using $(x+y)(x-y) = x^2 - y^2\,$, the product of the two denominators becomes
\begin{align*}
D_+ \, D_- &= q^2t^2\left( ((1+q)t - 2V\tilde{V})^2 -4R^2 \tilde{R}^2 \right) \\
&= q^2 t^2 \left( 4V^2 \tilde{V}^2 -4 R^2\tilde{R}^2 -4(1+q)tV\tilde{V} + t^2   \right)\\
&= q^2 t^2 \left(\frac{\lambda^2 \mu^2}{4} -4(t - \frac{\lambda^2}{4})(qt - \frac{\mu^2}{4}) -(1+q)t\lambda \mu + (1+q)^2t^2   \right) \\
&=q^2t^2\left(t\mu^2 +   qt \lambda^2 -(1+q)t \lambda \mu -4qt^2 + (1+q)^2t^2  \right)\\
&= q^2t^3 \left((\mu - \lambda)(\mu - q\lambda) + (1-q)^2t\right) \,,
\end{align*}
so that finally,
\begin{equation*}
\Re(S_+ - S_-) = \frac{2q^2(1-q)t^2R\tilde{R}}{q^2t^3 \left((\mu - \lambda)(\mu - q\lambda) + (1-q)^2t\right)} = \frac{2\pi^2 q(1-q) t \rho(\lambda,t) \tilde{\rho}(\mu,t)}{(1-q)^2t + (\mu - \lambda)(\mu - q\lambda)} \,,
\end{equation*}
and we obtain the announced formula (\ref{eq:solution_noise_only}):
\begin{equation*}
W(\mu, \lambda,t) = \frac{(1-q)\,t}{(1-q)^2\,t + (\mu - \lambda)\,(\mu - q\lambda)} \,.
\end{equation*}

\subsection{Study of the cubic equation}
\label{subsec:App_cubic}
We define the polynomial $P\left(X\right) := q\, X^3 - \left((1 + 6q + q^2)\, t + \mu^2\right) X + 4(1+q) \, t \, \mu\,$. One can first check that
\begin{align*}
P\left(2 \sqrt{t}\right) = - 2 \sqrt{t} \, \left((1+q)\, \sqrt{t} - \mu\right)^2 < 0 \,, \\
P\left(-2 \sqrt{t}\right) = 2 \sqrt{t} \, \left((1+q) \, \sqrt{t} + \mu\right)^2 > 0 \,.
\end{align*}
These inequalities can only become large in extreme cases like $q = 1$ or $t = 0\,$, because we always have $-2\sqrt{qt} \leq \mu \leq 2\sqrt{qt}\,$.
Since $\lim_{x \to \pm \infty} \, P(x) = \pm \infty \,$, we conclude that $P$ has three real roots: one in $]-\infty \, , -2\sqrt{t}[\,$, one in $[-2\sqrt{t}\, , 2 \sqrt{t}]\,$, and one in $]2\sqrt{t}\,,+\infty[\,$. Therefore, our maximisation problem always has a unique solution $\lambda_* \in [-2\sqrt{t}\,,2\sqrt{t}]\,$.

Recalling that $\mu = \mu(x,t)\,$, we want to prove that
$$
\lambda(qx + 1 -q,t) \leq \lambda_* \leq \lambda(qx,t) \,.
$$
In order to do so, we first demonstrate some properties of the quantile functions $\lambda(\cdot,t)$ and $\mu(\cdot,t)\,$:
\begin{itemize}
\item[$\bullet$] Since $\rho(\cdot, t)$ is even, one can easily check that $\lambda(1-x, t) = - \lambda(x,t)\,$.
\item[$\bullet$] We have the scaling inequality $\sqrt{q}\, \lambda(x,t) \leq \lambda(qx,t)$ for $q \in [0 \,, 1]\,$. This comes from comparing the integrals that define the quantile function (see (\ref{eq:quantile_function})):
\begin{align*}
\int_{\lambda(qx,t)}^{+\infty} \rho(\lambda,t) \, d\lambda &= qx \\
&= q\, \int_{\lambda(x,t)}^{+\infty} \frac{\sqrt{4t - \lambda^2}}{2\pi t} \, d\lambda \\
&= \sqrt{q} \, \int_{\sqrt{q}\, \lambda(x,t)}^{+\infty} \frac{\sqrt{4t - \lambda^2 \, / \, q}}{2\pi t} \, d\lambda \\
&= \int_{\sqrt{q} \, \lambda(x,t)}^{+\infty} \frac{\sqrt{4qt - \lambda^2}}{2\pi t} \, d\lambda \\
&\leq \int_{\sqrt{q} \, \lambda(x,t)}^{+\infty} \frac{\sqrt{4t - \lambda^2}}{2 \pi t} \, d\lambda \,.
\end{align*}
Since the integrand is the same in the initial and final integrals and is positive, we deduce that $\sqrt{q} \, \lambda(x,t) \leq \lambda(qx,t)\,$.
\item[$\bullet$] We have the relation $\mu(x,t) = \sqrt{q} \, \lambda(x,t)$ between the quantile functions associated with $\tilde{\rho}$ and with $\rho\,$. Once again, we use their definitions:
\begin{align*}
\int_{\mu(x,t)}^{+\infty} \frac{\sqrt{4qt - \mu^2}}{2 q \pi t} \, d\mu &= x \\
&= \int_{\lambda(x,t)}^{+\infty} \frac{\sqrt{4t - \lambda^2}}{2 \pi t} \, d\lambda \\
&= \int_{\sqrt{q} \, \lambda(x,t)}^{+\infty} \frac{\sqrt{4t - \lambda^2 \,/\, q}}{2 \sqrt{q} \pi t} \, d\lambda \\
&= \int_{\sqrt{q} \, \lambda(x,t)}^{+\infty} \frac{\sqrt{4qt - \lambda^2}}{2 q \pi t} \, d\lambda \,.
\end{align*} 
The integrands in the initial integral is the same as in the final one. Additionally, since we work inside the support of $\tilde{\rho}\,$, the integrand is strictly positive. Therefore, the integration bounds must be identical, i.e. $\mu(x,t) = \sqrt{q} \, \lambda(x,t)\,$.
\end{itemize}
We now consider the case $x \geq 1\,/\,2\,$, meaning that $\mu \geq 0$ and $\lambda(x,t) \geq 0\,$. The case $x < 1\,/\,2$ is symmetric. Combining the first two properties allows us to obtain
\begin{equation*}
\lambda(qx + 1 - q,t) = \lambda(1 - q\,(1-x),t) = - \lambda(q\,(1-x),t) \leq -\sqrt{q} \, \lambda(1-x,t) = \sqrt{q} \, \lambda(x,t) \,.
\end{equation*}
Since $\lambda(x,t) \geq 0$ and $\lambda(\cdot,t)$ is non-increasing,
\begin{equation*}
\lambda(qx + 1 - q,t) \leq \sqrt{q} \, \lambda(x,t) \leq \lambda(x,t) \leq \lambda(qx,t) \,,
\end{equation*}
or, using our third property,
\begin{equation*}
\lambda(qx+1-q,t) \leq \mu \leq \frac{\mu}{\sqrt{q}} \leq \lambda(qx,t) \,.
\end{equation*}
Finally, one can check that
\begin{equation*}
P(\mu) = (1-q) \mu \, \left((3+q)\,t - \mu^2 \right) \geq 3(1-q)^2 \mu t \geq 0 \,,
\end{equation*}
using the fact that $\mu^2 \leq 4qt\,$, and that
\begin{equation*}
P\left(\frac{\mu}{\sqrt{q}}\right) = - \frac{\mu t}{\sqrt{q}} \left(1 - \sqrt{q}\right)^4 \leq 0 \,.
\end{equation*}
Since both $\mu$ and $\mu \, / \, \sqrt{q}$ lie in $[-2\sqrt{t} \,, 2\sqrt{t}]\,$, we must have $\mu \leq \lambda_* \leq \mu \,/\, \sqrt{q}\,$, i.e.
\begin{equation*}
\lambda(qx+1-q,t) \leq \lambda_* \leq \lambda(qx,t) \,.
\end{equation*}

\subsection{Spike-spike overlap equation}
\label{subsec:App_spike_spike}
For this proof, we will need a spike variant of the reduction property (\ref{eq:property_R_bulk}), using the fact that $\lambda_1^N(t)$ is not complex but converges to $\lambda_1(t)\,$, which is distinct from the limiting bulk of radius $2\sqrt{t}\,$:
\begin{itemize}
\item[$\bullet$] \underline{Reduction property for the spike case:}

For any $z \in \mathbb{C} \setminus \mathbb{R}$ and for any $ p \geq 0\,$,
\begin{equation}
\tag	{R'}
\label{eq:property_R_spike}
\frac{1}{N} \, \sum_{k > 1}^N \, \frac{1}{|\lambda_1^N(t) - \lambda_k|^p} \, \left| \sum_{i = 1}^n \, \frac{\bk{i}{k} \bk{i}{1}}{z - \mu_i} \right|^2 = \mathcal{O}\left(\frac{1}{\sqrt{N}}\right) \,,
\end{equation}
\begin{equation}
\tag{R''}
\label{eq:property_R_spike_double}
\frac{1}{N^2} \, \sum_{\substack{j,k > 1 \\ k \neq j}}^N \, \frac{1}{|\lambda_1^N(t) - \lambda_j|^p \, |\lambda_1^N(t) - \lambda_k|^p} \,\left|\sum_{i=1}^n \, \frac{\bk{i}{j} \bk{i}{k}}{z - \mu_i}\right|^2 = \mathcal{O}\left(\frac{1}{\sqrt{N}}\right) \,.
\end{equation}

\noindent
\textit{Proof}

We prove the first equation (\ref{eq:property_R_spike}).
Let $\Sigma_p$ be the sum under investigation. For any $p \geq 0\,$, $\Sigma_p$ is of order 1, as the denominator $|\lambda_1^N(t) - \lambda_k|$ does not diverge in the scaling limit since $\lambda_1(t)$ is a spike. 

Next, we note the following inequality:
\begin{align*}
\frac{1}{|\lambda_1^N(t) - \lambda_k|^p} - \frac{1}{|\lambda_1^N(t) - \lambda_k|^p + \frac{1}{\sqrt{N}}} &= \frac{1}{\sqrt{N}\,|\lambda_1^N(t) - \lambda_k|^p \left(|\lambda_1^N(t) - \lambda_k|^p + \frac{1}{\sqrt{N}}\right)} \\
&\leq \frac{1}{\sqrt{N} \,|\lambda_1^N(t) - \lambda_k|^{2p}} \,,
\end{align*}
so that
\begin{align*}
\Sigma_p &\leq \frac{1}{N} \, \sum_{\substack{k > 1}}^N \frac{1}{|\lambda_1^N(t) - \lambda_k|^p + \frac{1}{\sqrt{N}}} \, \left|\sum_{i=1}^n \, \frac{\bk{i}{k} \bk{i}{1}}{z - \mu_i} \right|^2 + \frac{1}{\sqrt{N}} \, \Sigma_{2p} \\
&\leq \frac{1}{N} \, \sum_{\substack{k > 1}}^N \, \frac{1}{|\lambda_1^N(t) - \lambda_k|^p + \frac{1}{\sqrt{N}}} \, \left|\sum_{i=1}^n \, \frac{\bk{i}{k} \bk{i}{1}}{z - \mu_i} \right|^2 + \mathcal{O}\left(\frac{1}{\sqrt{N}}\right) \,.
\end{align*}
The idea here is to approximate $\lambda_1^N(t)$ by $\lambda_1^N(t) + i \, \varepsilon_N\,$, where $\varepsilon_N$ tends to zero at an optimal speed. Therefore, we can apply an upper bound similar to the one used in (\ref{eq:property_R_bulk}). Since $\frac{1}{|\lambda_1^N(t) - \lambda_j|^p + \frac{1}{\sqrt{N}}}$ is bounded from above by $\sqrt{N}\,$, we proceed similarly to the proof of (\ref{eq:property_R_bulk}). Denoting by $\tilde{\Sigma}_p$ the sum
$$
\frac{1}{N} \, \sum_{\substack{k > 1}}^N  \, \frac{1}{|\lambda_1^N(t) - \lambda_k|^p + \frac{1}{\sqrt{N}}} \, \left|\sum_{i=1}^n \, \frac{\bk{i}{k} \bk{i}{1}}{z - \mu_i} \right|^2 \,,
$$
we have
\begin{align*}
\tilde{\Sigma}_p &\leq \frac{1}{\sqrt{N}} \,  \sum_{\substack{k > 1}}^N \, \left|\sum_{i=1}^n \, \frac{\bk{i}{k} \bk{i}{1}}{z - \mu_i} \right|^2 \\
&\leq \frac{1}{\sqrt{N}} \, \sum_{\substack{i,l=1}}^n \, \frac{\bk{i}{1} \bk{l}{1}}{(z - \mu_i)(z^* - \mu_{l})} \, \sum_{k>1}^N \, \bk{i}{k} \bk{l}{k} \\
&\leq \frac{1}{\sqrt{N}} \, \sum_{\substack{i,l=1}}^n \, \frac{\bk{i}{1} \bk{l}{1}}{(z - \mu_i)(z^* - \mu_{l})} \, (\delta_{il} - \bk{i}{1} \bk{l}{1}) \\
&\leq \frac{1}{\sqrt{N}} \,  \left(\sum_{i=1}^n \, \frac{\bk{i}{1}^2}{|z - \mu_i|^2} - \left| \sum_{i=1}^n \, \frac{\bk{i}{1}}{z - \mu_i}\right|^2 \right) \\
&\leq \frac{1}{\Im(z)^2 \,\sqrt{N}} \,.
\end{align*}

Finally, we conclude with:
\begin{equation*}
\Sigma_p = \mathcal{O}\left(\frac{1}{\sqrt{N}}\right) \,.
\end{equation*}
One can verify that this is the best rate of convergence possible using this method. Specifically, using $(|\lambda_1^N(t) - \lambda_j|^p + N^{-\alpha})^{-1}$ results in a order of $\mathcal{O}(N^{-\alpha} \,+\, N^{\alpha - 1})\,$, which is optimal for $\alpha = 1\,/\,2\,$.

Finally, property (\ref{eq:property_R_spike_double}) is easily proved using the same method, by comparing $(|\lambda_1^N(t) - \lambda_j|^p \, |\lambda_1^N(t) - \lambda_k|^p)^{-1}$ and $(|\lambda_1^N(t) - \lambda_j|^p \, |\lambda_1^N(t) - \lambda_k|^p \,+\, N^{-1/2})^{-1}\,$.
\end{itemize}

\noindent
We first note that property (\ref{eq:property_R_spike_double}) allows us to extend the results from Appendix \ref{subsec:stieltjes_equation} to the spike case. Specifically,
\begin{equation*}
\frac{1}{N} \,\sum_{i=1}^n \, \sum_{j > 1}^N \, \frac{\bk{i}{j}^2}{(\lambda_1^N(t) - \lambda_j)(z - \mu_i)} \to S\left(\lambda_1(t),z, t\right) = \frac{q}{(\lambda_1(t) - t\,G(\lambda_1(t),t))\, (z-qt\,\tilde{G}(z,t)) - qt} \,.
\end{equation*} 
Indeed, (\ref{eq:property_R_spike_double}) adapts the original property (\ref{eq:property_R_bulk}), which was the only point for which we used the fact that $z \notin \mathbb{R}$ in the demonstration from Appendix \ref{subsec:stieltjes_equation} (one would get an order $1 \,/\, N^{3/2}$ variance for $dI_W$ instead of $1 \,/\, N^2$, but this is still summable so the almost sure convergence holds). One can also verify that (\ref{eq:property_R_spike}) and (\ref{eq:property_R_spike_double}) hold with $z \leftarrow \mu_1^N(t)\,$, which is the form we will use in this computation.

Using our Itô formula (\ref{eq:ito_squared_overlaps}) for the squared overlaps, we have
\begin{align}
\label{eq:ito_overlap_spike_spike}
d\left(\bk{1}{1}^2\right) &= \frac{1}{N} \, \sum_{\substack{k > 1}}^N \,  \frac{\bk{1}{k}^2 - \bk{1}{1}^2}{(\lambda_1^N(t) - \lambda_k)^2} \, dt + \frac{1}{N} \, \sum_{\substack{l > 1}}^n \, \frac{\bk{l}{1}^2 - \bk{1}{1}^2}{(\mu_1^N(t) - \mu_l)^2} \, dt \nonumber \\
&\quad + \frac{2}{N} \, \sum_{l > 1}^n \, \sum_{k > 1}^N \, \frac{\bk{1}{1}^2 \bk{l}{k}^2 + \bk{1}{k}^2 \bk{l}{1}^2 + 2 \bk{1}{1} \bk{l}{k} \bk{1}{k} \bk{l}{1}}{(\mu_1^N(t) - \mu_l)(\lambda_1^N(t) - \lambda_{k})} \, dt \nonumber \\
&\quad + \frac{2}{\sqrt{N}} \, \sum_{\substack{k>1}}^N \, \frac{dW_{1k}}{\lambda_1^N(t) - \lambda_k} \, \bk{1}{k} \bk{1}{1} + \frac{2}{\sqrt{N}} \, \sum_{\substack{l > 1}}^n \, \frac{d\tilde{W}_{1l}}{\mu_1^N(t) - \mu_l} \, \bk{l}{1} \bk{1}{1} \,.
\end{align}
We will demonstrate that this equation converges to (\ref{eq:f_equation_not_simplified}) without any rescaling, provided that all the overlaps involving an eigenvector in the bulk are of order $1\,/\,N$ (which makes sense considering the normalization constraints). This would confirm that $\bk{1}{1}^2$ remains of order 1 in the scaling limit. We begin by recognizing the following convergences:
\begin{itemize}
\item[$\bullet$] $\frac{1}{N} \, \sum\limits_{\substack{k>1}}^N \, \frac{\bk{1}{1}^2}{(\lambda_1^N(t) - \lambda_k)^2} \,\longrightarrow  \, \int_{\mathbb{R}} \frac{\rho(\lambda,t)}{(\lambda_1(t) - \lambda)^2}\, d\lambda \, f(t)\,$, where $\rho$ is the semicircular density of radius $2\sqrt{t}\,$, corresponding to the limiting density of all the $\lambda_k$ with $k \neq 1\,$.
\item[$\bullet$] $\frac{1}{N} \, \sum\limits_{l>1}^n \, \frac{\bk{1}{1}^2}{(\mu_1^N(t) - \mu_l)^2} \, \longrightarrow  \, \int_{\mathbb{R}} \frac{\tilde{\rho}(\mu,t)}{(\mu_1(t) - \mu)^2} \, d\mu \, f(t)\,$.
\item[$\bullet$] $\frac{1}{N} \, \sum\limits_{l> 1}^n \, \sum\limits_{k>1}^N \, \frac{\bk{1}{1}^2\bk{l}{k}^2}{(\mu_1^N(t) - \mu_l)(\lambda_1^N(t) - \lambda_{k})} \, \longrightarrow \, S\left(\lambda_1(t), \mu_1(t), t\right)\, f(t)\,$, from what we just proved.
\end{itemize}
All the other $dt$-terms vanish in the scaling limit. The sum
\begin{equation*}
\sum_{\substack{k>1}}^N \,  \frac{\bk{1}{k}^2}{(\lambda_1^N(t) - \lambda_k)^2}
\end{equation*}
remains of order 1 as $N$ increases because $\bk{1}{k}^2$ is of order $1\,/\,N$ and $\lambda_1^N(t)$ is asymptotically distinct from the bulk for $t < \lambda^2\,$. Since it is multiplied by $1\,/\,N$ in (\ref{eq:ito_overlap_spike_spike}), it converges to 0. The same reasoning applies to
\begin{equation*}
\frac{1}{N} \, \sum_{l > 1}^n \, \frac{\bk{l}{1}^2}{(\mu_1^N(t) - \mu_l)^2} = \mathcal{O}\left(\frac{1}{N}\right) \,,
\end{equation*}
and to
\begin{equation*}
\frac{4}{N} \, \sum_{l > 1}^n \, \sum_{k > 1}^N \, \frac{\bk{l}{k}\bk{1}{k}\bk{l}{1}}{(\mu_1^N(t) - \mu_l)(\lambda_1^N(t) - \lambda_{k})} = \mathcal{O}\left(\frac{1}{\sqrt{N}}\right) \,.
\end{equation*}
The argument is similar for the Brownian terms, but we must analyse their variances. Since they are centered, the variance of the first Brownian sum is
\begin{equation*}
\frac{4}{N}  \, \sum_{\substack{k,k'>1}}^N \,  \mathbb{E}\left[ \frac{dW_{1k} dW_{1k'}\,\bk{1}{k}\bk{1}{k'} \bk{1}{1}^2}{(\lambda_1^N(t) - \lambda_k)(\lambda_1^N(t) - \lambda_{k'})} \right] = \frac{4}{N}  \, \sum_{\substack{k>1}}^N \, \mathbb{E}\left[ \frac{\bk{1}{k}^2 \bk{1}{1}^2}{(\lambda_1^N(t) - \lambda_k)^2}\right]\,dt = \mathcal{O}\left(\frac{1}{N}\right) \,,
\end{equation*}
where we only used the fact that $\bk{1}{k}^2$ is of order $1\,/\,N$ and $\bk{1}{1}^2$ is bounded from above by 1. Therefore, the sum converges to 0 in $L^2\,$, and the same argument applies to the second Brownian sum.

Finally, the only terms that remain give us the deterministic limiting equation (\ref{eq:f_equation_not_simplified}):
\begin{equation}
f'(t) = \left[ - \int_{\mathbb{R}} \frac{\rho(\lambda,t)}{(\lambda_1(t) - \lambda)^2}\,d\lambda - q\, \int_{\mathbb{R}} \frac{\tilde{\rho}(\mu,t)}{(\mu_1(t) - \mu)^2}\,d\mu + 2 \, S\left(\lambda_1(t), \mu_1(t), t\right) \right] \, f(t) \,.
\end{equation}
Note that all these convergences are justified by the fact that $\lambda_1^N(t)$ and $\mu_1^N(t)$ asymptotically become the spikes $\lambda_1(t)$ and $\mu_1(t)\,$, which are distinct from their respective bulks. Therefore, the term $1\,/\,N \, \sum_k \, 1 \,/\, (\lambda_1^N(t) - \lambda_k)^2$ does not diverge in the scaling limit (similarly for $\mu_1^N(t)$).

\subsection{Spike-bulk Stieltjes transform equation}
\label{subsec:App_spike_bulk}
With what we proved in Appendix \ref{subsec:App_spike_spike}, we also have the convergence 
\begin{align*}
\frac{1}{N} \sum_{\substack{i = 1 \\ j > 1}}^{n,N} \frac{\bk{i}{j}^2}{(\lambda_1^N(t) - \lambda_j)^2(z - \mu_i)} &\longrightarrow - \partial_1 S(\lambda_1(t), z, t) \\
&= \frac{q(1 - t\partial_1 G(\lambda_1(t),t))(z - qt\tilde{G}(z,t))}{\left((\lambda_1(t) - tG(\lambda_1(t),t))(z-qt\tilde{G}(z,t)) - qt \right)^2} \,.
\end{align*}
Using the simplifications derived in \ref{subsec:spike_spike}, we obtain
\begin{equation}
\label{eq:limit_S_spike}
\frac{1}{N} \sum_{\substack{i = 1 \\ j > 1}}^{n,N} \frac{\bk{i}{j}^2}{(\lambda_1^N(t) - \lambda_j)(z - \mu_i)} \longrightarrow \frac{q}{\lambda (z - qt \tilde{G}(z,t)) - qt} \,,
\end{equation}
and
\begin{equation}
\label{eq:limit_dS_spike}
\frac{1}{N} \sum_{\substack{i=1 \\ j > 1}}^{n,N} \frac{\bk{i}{j}^2}{(\lambda_1^N(t) - \lambda_j)^2(z - \mu_i)} \longrightarrow \frac{q \lambda^2 (z - qt\tilde{G}(z,t))}{\left(\lambda^2 - t\right)\left(\lambda(z - qt \tilde{G}(z,t)) - qt \right)^2} \,.
\end{equation}

Now, we apply Itô's lemma to $S_{\lambda}^{(n)}\,$, without explicitly writing the terms that vanish as $N$ goes to infinity (see Appendix \ref{subsubsec:deriving_S_equation} for a detailed treatment of these terms):
\begin{equation}
\label{eq:ito_S_lambda}
dS_{\lambda}^{(n)} = \sum_{\substack{i = 1}}^n \, \frac{d(\bk{i}{1}^2)}{z - \mu_i} + \sum_{\substack{i = 1}}^n \, \frac{\bk{i}{1}^2}{(z - \mu_i)^2}\,d\mu_i + o(1) \,.
\end{equation}
Let us decompose the first sum using the Itô formula (\ref{eq:ito_squared_overlaps}),
\begin{equation*}
\sum_{\substack{i = 1}}^n\, \frac{d(\bk{i}{1}^2)}{z - \mu_i} = (I'_{\mu} + I'_{\lambda} + I'_{\lambda \mu})\,dt + dI'_W + dI'_{\tilde{W}} \,,
\end{equation*}
where:
\begin{itemize}
\item[$\bullet$] $I'_{\mu} = \frac{1}{N} \,\sum\limits_{\substack{i,l=1 \\ l \neq i}}^n \,\frac{\bk{l}{1}^2 - \bk{i}{1}^2}{(\mu_i - \mu_l)^2(z - \mu_i)}\,,$
\item[$\bullet$] $I'_{\lambda} = \frac{1}{N}\, \sum\limits_{\substack{i=1 \\ j > 1}}^{n,N} \,\frac{\bk{i}{j}^2 - \bk{i}{1}^2}{(\lambda_1^N(t) - \lambda_j)^2(z - \mu_i)}\,,$
\item[$\bullet$] $I'_{\lambda \mu} = \frac{2}{N}\, \sum\limits_{\substack{i,l = 1 \\ s \neq i }}^{n} \, \sum\limits_{j > 1}^N \frac{\bk{i}{1}^2 \bk{l}{j}^2 + \bk{i}{j}^2 \bk{l}{1}^2 + 2  \,\bk{i}{1} \bk{l}{j} \bk{i}{j}\bk{l}{1}}{(\mu_i - \mu_l)(\lambda_1^N(t) - \lambda_j)(z- \mu_i)}\,,$
\item[$\bullet$] $dI'_W = \frac{2}{\sqrt{N}} \,\sum\limits_{\substack{i = 1 \\ j > 1}}^{n,N} \,\frac{dW_{1j} \,\bk{i}{j} \bk{i}{1}}{(\lambda_1^N(t) - \lambda_j)(z - \mu_i)}\,,$
\item[$\bullet$] $dI'_{\tilde{W}} = \frac{2}{\sqrt{N}}\, \sum\limits_{\substack{i,l = 1 \\ l \neq i}}^n \,\frac{d\tilde{W}_{il}\, \bk{l}{1} \bk{i}{1}}{(\mu_i - \mu_l)(z - \mu_i)}\,.$
\end{itemize}
We explicit the limit of each of these terms in order to explain how equation (\ref{eq:S_lambda_equation}) arises. The manipulations are very similar to those in Appendix \ref{subsubsec:deriving_S_equation}:
\begin{enumerate}
\item $I'_{\mu}\,dt$ combined with the second sum in (\ref{eq:ito_S_lambda}) converges to $-q \,\tilde{G}(z,t) \, \partial_zS_{\lambda}\,dt\,$.
\item $I'_{\lambda}$ converges to $\frac{q\lambda^2(z - qt\tilde{G}(z,t))}{\left(\lambda^2 - t \right)\left(\lambda (z - qt\tilde{G}(z,t) -qt \right)^2} - \int \frac{\rho(\lambda,t)}{(\lambda_1(t) - \lambda)^2}\,d\lambda \,S_{\lambda}(z,t)\,$, using (\ref{eq:limit_dS_spike}). The integral simplifies to $1 \,/\, (\lambda^2 - t)\,$.
\item $I'_{\lambda \mu}$ converges to $\frac{2q}{\lambda (z - qt\tilde{G}(z,t)) - qt} \, S_{\lambda}(z,t)\,$, using the $(i,s)$ symmetry and equation (\ref{eq:limit_S_spike}), as well as property (\ref{eq:property_R_spike}) which shows that when we expand the sum in three terms, the third one vanishes.
\item $dI'_W$ converges to 0 in the scaling limit using (\ref{eq:property_R_spike}). This time the variance can be shown to be of order $1 \,/\, \sqrt{N}\,$, so almost sure convergence is not guaranteed.
\item $dI'_{\tilde{W}}$ also converges to 0. This can be shown using the same method we applied to $dI_W$ and $dI_{\tilde{W}}$ in Appendix \ref{subsubsec:deriving_S_equation}.
\end{enumerate}
Combining all these terms gives us equation (\ref{eq:S_lambda_equation}):
\begin{align*}
\partial_t S_{\lambda} &= -q \tilde{G}(z,t) \,\partial_z S_{\lambda} + \left(\frac{2q}{\lambda (z - qt \tilde{G}(z,t)) - qt} + \frac{1}{t - \lambda^2} \right) S_{\lambda} \\
&\quad + \frac{q \lambda^2 (z - qt\tilde{G}(z,t))}{\left(\lambda^2 - t\right)\left(\lambda (z - qt\tilde{G}(z,t)) - qt\right)^2} \,.
\end{align*}

\bibliographystyle{plain}
\bibliography{ref}

\begin{thebibliography}{10}

\bibitem{adler2013random}
Mark Adler, Pierre Van~Moerbeke, and Dong Wang.
\newblock Random matrix minor processes related to percolation theory.
\newblock {\em Random Matrices: Theory and Applications}, 2(04):1350008, 2013.

\bibitem{allez2014eigenvectors}
Romain Allez, Jo{\"e}l Bun, and Jean-Philippe Bouchaud.
\newblock The eigenvectors of gaussian matrices with an external source.
\newblock {\em arXiv preprint arXiv:1412.7108}, 2014.

\bibitem{anderson2010introduction}
Greg~W Anderson, Alice Guionnet, and Ofer Zeitouni.
\newblock {\em An introduction to random matrices}.
\newblock Number 118. Cambridge university press, 2010.

\bibitem{bai2010spectral}
Zhidong Bai and Jack~W Silverstein.
\newblock {\em Spectral analysis of large dimensional random matrices},
  volume~20.
\newblock Springer, 2010.

\bibitem{baryshnikov2001gues}
Yu~Baryshnikov.
\newblock Gues and queues.
\newblock {\em Probability Theory and Related Fields}, 119:256--274, 2001.

\bibitem{beenakker1994random}
CWJ Beenakker and B~Rejaei.
\newblock Random-matrix theory of parametric correlations in the spectra of
  disordered metals and chaotic billiards.
\newblock {\em Physica A: Statistical Mechanics and its Applications},
  203(1):61--90, 1994.

\bibitem{biane1997free}
Philippe Biane.
\newblock On the free convolution with a semi-circular distribution.
\newblock {\em Indiana University Mathematics Journal}, pages 705--718, 1997.

\bibitem{bolla2008noisy}
Marianna Bolla.
\newblock Noisy random graphs and their laplacians.
\newblock {\em Discrete Mathematics}, 308(18):4221--4230, 2008.

\bibitem{bourgade2016fixed}
Paul Bourgade, Laszlo Erd{\H{o}}s, Horng-Tzer Yau, and Jun Yin.
\newblock Fixed energy universality for generalized wigner matrices.
\newblock {\em Communications on Pure and Applied Mathematics},
  69(10):1815--1881, 2016.

\bibitem{bourgade2017eigenvector}
Paul Bourgade and H-T Yau.
\newblock The eigenvector moment flow and local quantum unique ergodicity.
\newblock {\em Communications in Mathematical Physics}, 350:231--278, 2017.

\bibitem{bun2018overlaps}
Jo{\"e}l Bun, Jean-Philippe Bouchaud, and Marc Potters.
\newblock Overlaps between eigenvectors of correlated random matrices.
\newblock {\em Physical Review E}, 98(5):052145, 2018.

\bibitem{cai2021asymptotic}
T~Tony Cai, Tiefeng Jiang, and Xiaoou Li.
\newblock Asymptotic analysis for extreme eigenvalues of principal minors of
  random matrices.
\newblock {\em The Annals of Applied Probability}, 31(6):2953--2990, 2021.

\bibitem{candes2005decoding}
Emmanuel~J Candes and Terence Tao.
\newblock Decoding by linear programming.
\newblock {\em IEEE transactions on information theory}, 51(12):4203--4215,
  2005.

\bibitem{denton2019eigenvectors}
Peter~B Denton, Stephen~J Parke, Terence Tao, and Xining Zhang.
\newblock Eigenvectors from eigenvalues.
\newblock {\em arXiv preprint ArXiv:1908.03795}, 2019.

\bibitem{deutsch1991quantum}
Josh~M Deutsch.
\newblock Quantum statistical mechanics in a closed system.
\newblock {\em Physical review a}, 43(4):2046, 1991.

\bibitem{drton2008moments}
Mathias Drton, H{\'e}l{\`e}ne Massam, and Ingram Olkin.
\newblock Moments of minors of wishart matrices.
\newblock 2008.

\bibitem{dyson1962brownian}
Freeman~J Dyson.
\newblock A brownian-motion model for the eigenvalues of a random matrix.
\newblock {\em Journal of Mathematical Physics}, 3(6):1191--1198, 1962.

\bibitem{feng2022principal}
Renjie Feng, Gang Tian, Dongyi Wei, and Dong Yao.
\newblock Principal minors of gaussian orthogonal ensemble.
\newblock {\em arXiv preprint arXiv:2205.05732}, 2022.

\bibitem{gatti2017applications}
Fabien Gatti, Benjamin Lasorne, Hans-Dieter Meyer, and Andr{\'e} Nauts.
\newblock {\em Applications of quantum dynamics in chemistry}, volume~98.
\newblock Springer, 2017.

\bibitem{gorin2017interacting}
Vadim Gorin and Mykhaylo Shkolnikov.
\newblock Interacting particle systems at the edge of multilevel dyson brownian
  motions.
\newblock {\em Advances in Mathematics}, 304:90--130, 2017.

\bibitem{grabiner1999brownian}
David~J Grabiner.
\newblock Brownian motion in a weyl chamber, non-colliding particles, and
  random matrices.
\newblock In {\em Annales de l'IHP Probabilit{\'e}s et statistiques},
  volume~35, pages 177--204, 1999.

\bibitem{granziol2022learning}
Diego Granziol, Stefan Zohren, and Stephen Roberts.
\newblock Learning rates as a function of batch size: A random matrix theory
  approach to neural network training.
\newblock {\em Journal of Machine Learning Research}, 23(173):1--65, 2022.

\bibitem{haemers1995interlacing}
Willem~H Haemers.
\newblock Interlacing eigenvalues and graphs.
\newblock {\em Linear Algebra and its applications}, 226:593--616, 1995.

\bibitem{hu2023extreme}
Jianwei Hu, Seydou Keita, and Kang Fu.
\newblock Extreme eigenvalues of principal minors of random matrices with
  moment conditions.
\newblock {\em Journal of the Korean Statistical Society}, 52(3):715--735,
  2023.

\bibitem{huang2022eigenvalues}
Jiaoyang Huang.
\newblock Eigenvalues for the minors of wigner matrices.
\newblock In {\em Annales de l'Institut Henri Poincare (B) Probabilites et
  statistiques}, volume~58, pages 2201--2215. Institut Henri Poincar{\'e},
  2022.

\bibitem{johansson2006eigenvalues}
Kurt Johansson and Eric Nordenstam.
\newblock Eigenvalues of gue minors.
\newblock 2006.

\bibitem{ledoit2011eigenvectors}
Olivier Ledoit and Sandrine P{\'e}ch{\'e}.
\newblock Eigenvectors of some large sample covariance matrix ensembles.
\newblock {\em Probability Theory and Related Fields}, 151(1):233--264, 2011.

\bibitem{li2015central}
Lingyun Li, Matthew Reed, and Alexander Soshnikov.
\newblock Central limit theorem for linear eigenvalue statistics for
  submatrices of wigner random matrices.
\newblock {\em arXiv preprint arXiv:1504.05933}, 2015.

\bibitem{najnudel2021eigenvector}
Joseph Najnudel.
\newblock Eigenvector convergence for minors of unitarily invariant infinite
  random matrices.
\newblock {\em International Mathematics Research Notices}, 2021(8):6293--6314,
  2021.

\bibitem{najnudel2021bead}
Joseph Najnudel and B{\'a}lint Vir{\'a}g.
\newblock The bead process for beta ensembles.
\newblock {\em Probability Theory and Related Fields}, 179(3):589--647, 2021.

\bibitem{o2016eigenvectors}
Sean O'Rourke, Van Vu, and Ke~Wang.
\newblock Eigenvectors of random matrices: a survey.
\newblock {\em Journal of Combinatorial Theory, Series A}, 144:361--442, 2016.

\bibitem{pacco2023overlaps}
Alessandro Pacco and Valentina Ros.
\newblock Overlaps between eigenvectors of spiked, correlated random matrices:
  From matrix principal component analysis to random gaussian landscapes.
\newblock {\em Physical Review E}, 108(2):024145, 2023.

\bibitem{potters2020first}
Marc Potters and Jean-Philippe Bouchaud.
\newblock {\em A first course in random matrix theory: for physicists,
  engineers and data scientists}.
\newblock Cambridge University Press, 2020.

\bibitem{tao2012topics}
Terence Tao.
\newblock {\em Topics in random matrix theory}, volume 132.
\newblock American Mathematical Soc., 2012.

\bibitem{tao2012random}
Terence Tao and Van Vu.
\newblock Random matrices: universal properties of eigenvectors.
\newblock {\em Random Matrices: Theory and Applications}, 1(01):1150001, 2012.

\bibitem{voiculescu1986addition}
Dan Voiculescu.
\newblock Addition of certain non-commuting random variables.
\newblock {\em Journal of functional analysis}, 66(3):323--346, 1986.

\end{thebibliography}

\end{document}